\documentclass[10pt]{article}
\usepackage{latexsym}
\usepackage{epic}
\usepackage{eepic}
\usepackage{color}
\usepackage{amsmath,multirow,amsfonts}
\usepackage{tikz-cd}
\usetikzlibrary{matrix}
\usepackage{array}
\usepackage[english]{babel}
\usepackage[latin1]{inputenc}
\usepackage{times}
\usepackage[T1]{fontenc}
\usepackage{amssymb}
\usepackage{graphics}
\usepackage{graphicx}
\usepackage{epsfig}
\usepackage{cite}
\usepackage{caption}
\usepackage{tikz}
\usetikzlibrary{arrows}
\tikzstyle{block}=[draw opacity=0.7,line width=1.4cm]
\begin{document}
\newtheorem{def1}{Definition}[section]
\newtheorem{lem}{Lemma}[section]
\newtheorem{exa}{Example}[section]
\newtheorem{thm}{Theorem}[section]
\newtheorem{pro}{Proposition}[section]
\newtheorem{cor}{Corollary}[section]
\newtheorem{rem}{Remark}[section]
\newtheorem{exam}{Example}[section]
\title{On categories of spaces with $L$-fuzzy partitions, $L$-fuzzy closure system spaces and coalgebras (dialgebras)}
\author{Abha Tripathi\thanks{tripathiabha29@gmail.com} and S.P. Tiwari\thanks{sptiwarimaths@gmail.com}\\
Department of Mathematics \& Computing\\ Indian Institute of Technology (ISM)\\
Dhanbad-826004, India}
\date{}
\maketitle
\begin{abstract}
In this contribution, we aim to introduce and study $L$-fuzzy partition spaces and $L$-fuzzy closure system spaces in a categorical framework. Further, we present the concepts of coalgebras and dialgebras corresponding to {a direct upper $F$-transform under certain conditions} and show the functorial relationship between the category of spaces with $L$-fuzzy partition and the category of coalgebras (dialgebras). Moreover, we show that the categories of coalgebras and dialgebras are isomorphic and introduce a pair of adjoint functors between the coalgebras and dialgebras.
\end{abstract}
\section{Introduction} 
The concept of fuzzy transform ($F$-transform) was firstly introduced by Perfilieva \cite{per}, a theory that attracted the interest of {many researchers}. It has now been significantly developed and opened a new page in the theory of semi-linear spaces. The main idea of the $F$-transform is to factorize (or fuzzify) the precise values of independent variables by a closeness relation, and precise values of dependent variables are averaged to an approximate value. The theory of $F$-transform has already been elaborated and extended from real-valued to lattice-valued functions (cf., \cite{{per},{irin1}}), from fuzzy sets to parametrized fuzzy sets \cite{st} and from the single variable to the two (or more variables) (cf., \cite{{ma},{mar}, {mar1}, {step}}). Recently, several researchers have initiated the study of $F$-transforms based on an arbitrary $L$-fuzzy partition of an arbitrary universe (cf., \cite{kh1,jir,mo, mocko, anan,spt1,rus,spt}), where $L$ is a complete residuated lattice. Among these studies, the relationships between $F$-transforms and {semimodule homomorphisms} were investigated in \cite{jir}; a categorical study of $L$-fuzzy partitions of an arbitrary universe was presented in \cite{mo}; while, the relationships between $F$-transforms and similarity { relations were} established in \cite{mocko}. Further, in \cite{anan}, an interesting relationship among $F$-transforms, $L$-fuzzy topologies/co-topologies and $L$-fuzzy approximation operators (which are concepts used in the study of an operator-oriented view of fuzzy rough set theory) was established, while in \cite{spt1}, the relationship between fuzzy pretopological spaces and spaces with $L$-fuzzy partition was established. Also, in a different direction, a generalization of $F$-transforms was presented in \cite{rus} by considering the so-called $Q$-module transforms, where $Q$ stands for an unital quantale, while $F$-transforms based on a generalized residuated lattice {were studied} in \cite{spt}. Further, classes of $F$-transforms taking into account the well-known classes of implicators, namely {$R-,S-, QL-$implicators} were discussed in \cite{tri}. The various studies carried out in the line of applications of $F$-transforms, e.g., trend-cycle estimation \cite{holc}, { data compression \cite{hut}}, numerical solution of partial differential equations \cite{kh}, scheduling \cite{li}, time series \cite{vil}, data analysis \cite{no}, denoising \cite{Ir}, face recognition \cite{roh}, neural network { approaches \cite{ste} and trading} \cite{to}.\\\\
The concept of category theory introduced by Eilenberg and Mac Lane \cite{eil} is well-known. Further, many researchers developed this theory in \cite{fre, gro,law,law1,pi}. In the framework of $F$-transforms, the category ${\bf SpaceFP}$ was introduced by Mo\v{c}ko\v{r} \cite{mo} as a generalization of categories of sets with fuzzy partitions defined by lattice valued fuzzy equivalences (or, equivalently, similarity relations). The fuzzy partition defined by lattice-valued fuzzy equivalences gives a bijective correspondence between fuzzy equivalences and fuzzy partitions. On the other hand, a fuzzy partition is defined as a system of fuzzy sets uniquely defined by the given fuzzy equivalence; however, it is, in many cases, a disadvantage. One of the examples where fuzzy partitions defined by fuzzy equivalence relation cannot be used in the theory of $L$-valued $F$-transform. The definition of fuzzy sets as elements of a fuzzy partition resulting from the given fuzzy equivalence relation does not allow these fuzzy sets to be independently modified, which is necessary, for example, to achieve the required accuracy of inverse $F$-transform. This is why fuzzy partitions defined by fuzzy relations were not chosen as the ground structure for the category of $F$-transform. In \cite{mock,jiri,jir,mocko,mockor}, the category ${\bf SpaceFP}$ of spaces with fuzzy partitions and some properties has been investigated. Also, it has been shown that ${\bf SpaceFP}$ is a topological category, and the categories of upper (lower) transformation systems, satisfying simplified axioms, are isomorphic to the category ${\bf SpaceFP}$. Moreover, a functorial relationship exists between a subcategory of ${\bf SpaceFP}$ and categories of Kuratowski closure and interior operators, and a category of fuzzy preorders, respectively. In another direction, categories of fuzzy topologies have also been introduced and studied in \cite{chen1,lai,qi,rama}, while the concepts of $L$-fuzzy closure system spaces and $L$-fuzzy closure spaces from the categorical point of view have been studied in \cite{fan}.\\\\
{Coalgebra is is a well-known abstract theory that provides a uniform framework for various transition systems arising as a relatively recent theory within or closely connected to category theory. Further, coalgebra has shown to be useful in describing mathematical structures such as automata, processes, labeled transition systems, probabilistic transition systems, modal logic, object-oriented design and component-based software engineering (cf., \cite{de,do,rutt, soko,sokol}). Moreover, a framework for coalgebraic semantics for quantum systems is presented in \cite{liu} and it is shown that the notion of bisimulation between coalgebras is straightforward and intuitive when it is induced by the coalgebraic view \cite{feng,larsen,hag}. For the study of fuzzy mathematical structures, coalgebra for fuzzy transition systems was studied in \cite{wu}, and coalgebra and dialgebra (a generalization of algebra and coalgebra, cf., \cite{poll}) for fuzzy automata have been studied in \cite{chen,sinha}. Considering that fuzzy automaton and fuzzy transition systems are fuzzy relational structures, it will be interesting to use the concepts of coalgebra and dialgebra to enrich the theory of $F$-transform. Interestingly, we found that the direct upper $F$-transform determines coalgebra and dialgebra under certain conditions. Also, we study the theory of $F$-transform with the concept of $L$-fuzzy closure system. It is to be pointed out here that the categorical studies carried out herein are based on a slightly different morphism having an exciting relationship with some of the familiar categories. Specifically, in the presented work, 
\begin{itemize}
\item we introduce the category of spaces with $L$-fuzzy partitions and show the existence of product in this category.
\item we establish a functorial relationship of the category $\textbf{FPS}$ (category of spaces with $L$-fuzzy partitions) with different categories such as ${\bf FAS}$ (category of $L$-fuzzy approximation spaces), $\textbf{FCSS}$ (category of $L$-fuzzy closure system spaces) and ${\bf FCS}$ (category of $L$-fuzzy closure spaces);
\item different functors on the category of sets lead us to show that the category ${\bf COA}$ (category of $T_1$-coalgebras) is isomorphic to the category ${\bf DIA}$ (category of $(T_2,T_3)$-dialgebras), where $T_1,T_2,T_3$ are functors; and
\item we establish the adjunction between the categories of coalgebras.
\end{itemize}}
The structure of the paper is as follows. 
Section 2 contains elementary knowledge about the content of the paper. In Section 3, we introduce the category $\textbf{FPS}$ and some basic properties of $F$-transform. Section 4 contains some well-known categories and the existence of functors among these categories. The concept of coalgebra and dialgebra corresponding to $L$-fuzzy partition is introduced. We also show a relationship between coalgebra and dialgebra.
\section{Preliminaries}
This section is divided into two subsections. We recall some basic notions and notations related to the category theory in the first subsection, while that of complete residuated lattices, $L$-fuzzy sets, $L$-fuzzy relations, and some other basic notions are recalled in the second subsection.
\subsection{Category theory}
{In this subsection, we recall the ideas associated with categories, $F$-coalgebras, $(F, G)$-dialgebras, which are used in the main text. For details, we refer \cite{ad,alt,ar,her,lan}. Throughout, for a category ${\bf C}$, $|{\bf C}|$ denotes the class of objects of ${\bf C}$; while its morphisms are written as ${\bf C}$-morphisms. For $C,D\in |{\bf C}|$, homset ${\bf C}(C,D)$ denotes the set of all morphisms from $C$ to $D$. We begin with the following.}
\begin{def1} 
An object $Z\in|{\bf C}|$ is called 
\begin{enumerate}
    \item[(i)] an {\bf initial object} if for each $C\in|{\bf C}|$ there is exactly one ${\bf C}$-morphism from $Z$ to $C$;
    \item[(ii)] a {\bf terminal object} if for each $C\in |{\bf C}|$ there is exactly one ${\bf C}$-morphism from $C$ to $Z$; and 
    \item[(iii)] a {\bf zero object} if it is both initial and terminal object. 
\end{enumerate}
\end{def1}
For $C,D\in|{\bf C}|$ and zero object $Z\in|{\bf C}|$, we use $0_{C,D}$ for the unique ${\bf C}$-morphism $C \rightarrow Z \rightarrow D$.
\begin{def1}
A {\bf semiadditive} category is a category ${\bf C}$, where for $B,C\in {\bf C}$, each homset ${\bf C}(B,C)$ is equipped with the structure of a commutative monoid with operation $`+$' such that for all $\phi : A \rightarrow B, \psi, \eta : B \rightarrow C$ and $ \zeta: C \rightarrow D$
\begin{itemize}
\item[(i)] $(\psi + \eta) \circ \phi = (\psi \circ \phi ) + (\eta \circ \phi )$; and
\item[(ii)] $\zeta \circ (\psi + \eta) = (\zeta \circ \psi) + (\zeta \circ \eta)$.
\end{itemize}
\end{def1}
\begin{def1}
For a functor $F : {\bf SET} \rightarrow {\bf SET}$, a $F$-{\bf coalgebra} is a pair $(X, \alpha)$, where $X \in |{\bf SET}|$ and $\alpha:
X \rightarrow F(X)$ is a structure function of $F$-coalgebra.
\end{def1}
\begin{def1}
For two functors $F, G : {\bf C} \rightarrow {\bf D}$, a $(F, G)$-{\bf dialgebra} is a pair $(X, \beta)$ , where $X \in |{\bf SET}|$ and $\beta : F(X) \rightarrow G(X)$ is a structure function of $(F, G)$-dialgebra.
\end{def1}
\begin{def1}\label{adj0}
Let ${\bf C}$, ${\bf D}$ be categories and $F:{\bf C}\rightarrow{\bf D},G:{\bf D}\rightarrow {\bf C}$ be functors. Then $F$ is {\bf left adjoint} to $G$ and $G$ is {\bf right adjoint} to $F$ if there exists a natural transformation $\Phi:id_{\bf C}\rightarrow G\circ F$ such that for every $C\in|{\bf C}|$, $D\in|{\bf D}|$ and ${\bf C}$-morphism $f:C\rightarrow G(D)$, there exists a unique ${\bf D}$-morphism $g:F(C)\rightarrow D$ such that the diagram in Figure \ref{fig,0} commutes.
\end{def1}
\begin{figure}
\[\begin{tikzcd}[row sep=12ex, column sep=12ex] 
   C\arrow{r}{\Phi_C} \arrow[swap]{dr}{f} & G(F(C))\arrow{d}{G(g)} &F(C)\arrow{d}{g}\\
     & G(D)&D
  \end{tikzcd}\]
 \caption{Diagram for Definition \ref{adj0}}
\label{fig,0}
  \end{figure}
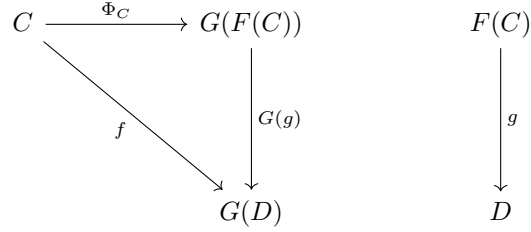
\subsection{Residuated lattice}
In this subsection, we recall some concepts related to complete residuated lattices and $L$-fuzzy relations, which are used in main text. For details, we refer \cite{abd,bel,belo,gog,kl,kli,ou}. We begin with the following.
\begin{def1}
A {\bf complete residuated lattice} is an algebra $(L,\wedge,\vee,\otimes, \rightarrow,0,1)$ such that 
\item(i) $(L,\wedge,\vee,0,1)$ is a complete lattice with the least element $0$ and the greatest element $1$;
\item(ii) $(L,\otimes,1)$ is a commutative monoid; and
\item (iii) $\forall a,b,c \in$ $L$,
\begin{center}
 $a\otimes b \leq c$ iff  $a\leq b\rightarrow c,$
 \end{center}
i.e., ($ \rightarrow ,\otimes$) is an adjoint pair on $L$.
\end{def1}
\begin{pro}
 Let $(L,\wedge,\vee,\otimes,\rightarrow,0,1)$ be a complete residuated lattice. Then for all $a,b,c\in L$, $\lbrace a_j,j\in J\rbrace\subseteq L$, the following properties hold:
 \begin{enumerate}
 \item[(i)] $a\otimes 0=0,\, a\otimes 1=a$;
 \item[(ii)] $a\rightarrow 1=1,\,  1\rightarrow a=a,a\rightarrow a=1$;
 \item[(iii)] If $a\leq b$ then $c\rightarrow a\leq c\rightarrow b,a\rightarrow c\geq b\rightarrow c$;
 \item[(iv)] $a\otimes(b\rightarrow c)=a\rightarrow (b\rightarrow c)$,  $a\otimes(b\rightarrow c)\leq (a\rightarrow b)\rightarrow c$;
 \item[(v)] $(a\rightarrow b)\otimes c\leq a\rightarrow b\otimes c$, $a\otimes b\rightarrow a\otimes c\geq b\rightarrow c$;
  \item[(vi)]$(a\rightarrow b)\rightarrow (c\rightarrow b)\geq c\rightarrow a$, $(a\rightarrow b)\rightarrow (a\rightarrow c)\geq b\rightarrow c$;
 \item[(vii)] $\bigvee\limits_{j\in J}a_j\otimes b=\bigvee\limits_{j\in J}(a_j\otimes b)$, $\bigwedge\limits_{j\in J}a_j\otimes b\leq\bigwedge\limits_{j\in J}(a_j\otimes b)$;
 \item[(viii)]  $a\rightarrow\bigwedge\limits_{j\in J}b_j=\bigwedge\limits_{j\in J}(a\rightarrow b_j)$, $\bigvee\limits_{j\in J}a_j\rightarrow b=\bigwedge\limits_{j\in J}(a_j\rightarrow b)$; and
  \item[(ix)] $\bigvee\limits_{j\in J}(a\rightarrow b_j)\leq a\rightarrow \bigvee\limits_{j\in J} b_j$.
 \end{enumerate}
 \end{pro}
 Further, if for all $a,b\in L,\,a\neq 0$ and $b\neq 0\Rightarrow a\otimes b\neq 0$, then $(L,\wedge,\vee,\otimes,\rightarrow,0,1)$ is called a {\bf complete residuated lattice without zero divisors}. Throughout this paper, we assume that $L\equiv(L,\wedge,\vee,\otimes,\rightarrow,0,1)$ is a complete residuated lattice without zero divisors and $L$-fuzzy sets are in the sense of \cite{gog}, i.e., an $L$-fuzzy set A of $X \in |{\bf SET}|,A:X\rightarrow L$. Further, for $X \in |{\bf SET}|, L^X$ denotes the collection of all $L$-fuzzy sets in $X$ and $J$ denotes an indexed set. Also, for all $a\in L,\textbf{{a}}\in L^X$ such that $\textbf{{a}}(x)=a,\,\forall\,x\in X$ denotes constant $L$-fuzzy set. For all $A,B\in L^X$, $A\leq B$ if $A(x)\leq B(x),\,\forall\,x\in X$. \\\\
The {\bf core} of an $L$-fuzzy set $A$ is defined as a crisp set
\begin{eqnarray*}
core(A)=\lbrace x\in X,A(x)=1\rbrace.
\end{eqnarray*}
If $core(A)\neq \emptyset$, then $A$ is called an {\bf normal $L$-fuzzy set}. Further, for $A\subseteq X\in|{\bf SET}|$, the characteristic function of $A$ is a function $1_A:X\rightarrow \{0,1\}$ such that
\begin{eqnarray*} 
1_A(x)=\begin{cases}
1 &\text{ if } x\in A,\\
0&\text{ otherwise}.
\end{cases} 
\end{eqnarray*}
Now, let $f : X \rightarrow Y$ be a $|{\bf SET}|$-morphism. Then according to Zadeh's extension principle, $f$ can be extended to the operators $\overrightarrow{f} : L^X \rightarrow L^Y$ and $\overleftarrow{f} : L^Y\rightarrow L^X$ such that for all $A \in L^X, B \in L^Y,y\in Y$
$$\overrightarrow{f}(A)(y)=\bigvee\limits_{x\in X,f(x)=y}A(x),\,\,\overleftarrow{f}(B)=B\circ f.$$
Also, for all $ X\in|{\bf SET}|,A,B\in L^X, x\in X$ and $\lbrace A_i:i\in J\rbrace\subseteq L^X$,
\begin{enumerate}
\item[(i)] $(A\otimes B)(x)=A(x)\otimes B(x)$;
\item[(ii)] $(A\rightarrow B)(x)=A(x)\rightarrow B(x)$;
 \item[(iii)] $(\bigwedge\limits_ {i\in J}A_i)(x)=\bigwedge\limits_ {i\in J}A_i(x)$; and
 \item[(iv)] $(\bigvee\limits_ {i\in J}A_i)(x)=\bigvee\limits_ {i\in J}A_i(x)$.
 \end{enumerate}
{In the following, we recall the concepts related to $L$-fuzzy relations.}
\begin{def1}\label{rela}
An {\bf $L$-fuzzy relation} on $X$ is a function $R:X\times X\rightarrow L$.
\end{def1}
\begin{def1}\label{app}
Let $R$ be an $L$-fuzzy relation on $X$. A pair $(X,R)$ is called an {\bf $L$-fuzzy approximation space}.
\end{def1}
\begin{def1}\label{lower}
Let $(X,R)$ be an $L$-fuzzy approximation space. Then the $L$-fuzzy upper approximation operator is a function $\overline{R}:L^X\rightarrow L^X$ such that
$$\overline{R}(f)(x)=\bigvee\limits_{y\in X}R(x,y)\otimes f(y),\,\forall\,x\in {X},f\in L^{X}.$$
\end{def1}
{We close this subsection, by recalling the concept of $L$-fuzzy closure operator.}
\begin{def1}\label{closure}
A function $c:L^X\rightarrow L^X$ is called an $L$-fuzzy closure operator if it satisfies:
\begin{itemize}
\item[(i)] $c(1_X)=1_X$;
\item[(ii)] $c(f)\geq f$;
\item[(iii)] $c(f\vee g)=c(f)\vee c(g)$; and
\item[(iv)] $c(c(f))= c(f)$.\\\\
The pair $(X,c)$ is called an {\bf $L$-fuzzy closure space}. 
\end{itemize}
\end{def1}
An $L$-fuzzy closure space $(X,c)$ is called {\bf strong}, if 
\[c(\textbf{a} \otimes f)\geq  \textbf{a}\otimes c(f) ,\,\forall\,\textbf{a},f\in L^X.\]
\section{Categorical view of $L$-fuzzy partitions and $L$-fuzzy relations}
In this section, we introduce and study the concepts of $L$-fuzzy partitions and $L$-fuzzy relations from the categorical point of view by using the morphisms in a slightly different way, which are generalization of the morphisms given in \cite{lai,mo}. Now, we recall the following from \cite{anan}.
\begin{def1}\label{FP}
A collection $\mathcal{P}$ of normal $L$-fuzzy sets $\lbrace A_{j}:j\in J\rbrace$ is called an {\bf $L$-fuzzy partition} of $X$ if the corresponding collection of ordinary sets $\lbrace core(A_{j}):j\in J\rbrace$ is partition of $X$. The pair $(X,\mathcal{P})$ is called a {\bf space with $L$-fuzzy partition}.
\end{def1}
Let $\mathcal{P} = \lbrace A_{j}:j\in J\rbrace$ be an $L$-fuzzy partition of $X$. Then it can be associated by the following onto index function $\xi:X\rightarrow J$, 
$$\xi(x)=j\Leftrightarrow x\in core(A_{j}).$$
\begin{rem} {(i) Let $A:X\rightarrow L$ and $B:Y\rightarrow L$ be two normal $L$-fuzzy sets in $X$ and $Y$, respectively, i.e., $core(A)\neq\emptyset$ and $core(B)\neq\emptyset$. We define an $L$-fuzzy set $A\times B:X\times Y\rightarrow L$ in $X\times Y$ such that $(A\times B)(x,y)=A(x)\wedge B(y),\,\forall\,x\in X,y\in Y.$ Because of nonempty core of $A$ and $B$, there exists $(x,y)\in X\times Y$ such that $(A\times B)(x,y)=A(x)\wedge B(y)=1$, or that $core(A\times B)\neq \emptyset$, showing that $A\times B$ is also a normal $L$-fuzzy set.}\\\\
(ii) Let $\lbrace core(A_{j}):j\in J_1\rbrace,\lbrace core(B_{j^{'}}):j^{'}\in J_2\rbrace$ and $\mathcal{P}_X=\{A_j:j\in J_1\},\mathcal{P}_Y=\{B_{j^{'}}:j^{'}\in J_2\}$ be partitions and $L$-fuzzy partitions of $X,Y$, respectively. Then $\lbrace core((A\times B)_{(j,j^{'})})=core(A_{j}\times B_{j^{'}}):(j,j^{'})\in J_1\times J_2,A_j\in\mathcal{P}_X,B_{j^{'}}\in\mathcal{P}_Y\rbrace$ is a partition of $X\times Y$ and $\mathcal{P}_{X\times Y}=\{(A\times B)_{(j,j^{'})}=A_j\times B_{j^{'}}:(j,j^{'})\in J_1\times J_2,A_j\in \mathcal{P}_X,B_{j^{'}}\in \mathcal{P}_Y\}$ is an $L$-fuzzy partition of $X\times Y$.
\end{rem}
\begin{exa}
Let $X=\mathbb{N}$ be set of natural numbers, $Y=\mathbb{Z}$ be set of integers and $L=[0,1]$. Then $ \mathcal{P}_X=\{A_1,A_2\},\mathcal{P}_Y=\{B_1,B_2\}$ are $L$-fuzzy partition of $X,Y$, respectively and for all $n\in\mathbb{N},m\in\mathbb{Z}$, $A_1,A_2:X\rightarrow L,B_1,B_2:Y\rightarrow L$ are defined as 
\begin{eqnarray*} 
A_1(n)=\begin{cases}
1 &{n~ is~ even},\\
0.4&{ n~ is~ odd},
\end{cases} 
A_2(n)=\begin{cases}
0.2 &{if~n ~is ~even},\\
1&{if ~n ~is~ odd},
\end{cases} 
\end{eqnarray*}
\begin{eqnarray*} 
B_1(m)=\begin{cases}
1 &{m ~is~ even},\\
0.4&{ m ~is ~odd},
\end{cases} 
B_2(m)=\begin{cases}
0.4 &{ m ~is ~even},\\
1&{~ m ~is~ odd}.
\end{cases} 
\end{eqnarray*}
Now, let $(A_j\times B_k)(x,y)=A_j(x)\wedge B_k(y),\,\forall\,x\in X,y\in Y,A_j\in\mathcal{P}_X,B_k\in\mathcal{P}_Y$. Then
\begin{eqnarray*}
core(A_1\times B_1)&=&\{(n,m)\in \mathbb{N}\times\mathbb{Z}:n\in\mathbb{N} ,m\in\mathbb{Z}~ are ~even\};\\
core(A_1\times B_2)&=&\{(n,m)\in \mathbb{N}\times\mathbb{Z}:n\in\mathbb{N} ~is ~even ~and~ m\in\mathbb{Z} ~is ~odd\};\\
 core(A_2\times B_1)&=&\{(n,m)\in \mathbb{N}\times\mathbb{Z}:n\in\mathbb{N}~ is~ odd ~and ~m\in\mathbb{Z} ~is ~even\};\,and\\
 core(A_2\times B_2)&=&\{(n,m)\in \mathbb{N}\times\mathbb{Z}:n\in\mathbb{N} ~is~ odd ~and ~m\in\mathbb{N}~is ~odd\}.
\end{eqnarray*}
Thus $\{core(A_j\times B_k):j=1,2,k=1,2\}$ form a partition of $\mathbb{N}\times \mathbb{Z}$ and $\mathcal{P}_{\mathbb{N}\times\mathbb{Z}}=\{A_j\times B_k:j=1,2,k=1,2\}$ is an $L$-fuzzy partition of $\mathbb{N}\times \mathbb{Z}$.
\end{exa}
{Chiefly inspired from the works in \cite{xing} regarding category of fuzzy automata, in the following, we introduce a mrophism betweeen $L$-fuzzy partitions for the categorical study.}
\begin{def1} \label{par}
Let $(X,\mathcal{P}_X)$ and $(Y,\mathcal{P}_Y)$ be two spaces with $L$-fuzzy partitions, where $\mathcal{P}_X= \{A_j:j\in J_1\}$ and $\mathcal{P}_Y= \{B_{j^{'}}:{j^{'}}\in J_2\}$. Then an {\bf FP-map} from $(X,\mathcal{P}_X)$ to $(Y,\mathcal{P}_Y)$ is a triple $(\phi,\psi,\mathcal{W})$ such that
\begin{itemize}
\item[(i)] $\phi:X\rightarrow Y, \,\psi:J_1\rightarrow J_2$ are functions;
\item[(ii)] $\mathcal{W}\subseteq \mathcal{P}_X\times \mathcal{P}_Y$ is a relation such that $\textbf{dom}(\mathcal{W})=\mathcal{P}_X$; and
\item[(iii)] $\exists\,l\in L\setminus 0$ such that for all $x\in X,j\in J_1,\,(A_j,B_{\psi(j)})\in \mathcal{W}$, 
$$l\rightarrow B_{\psi(j)}(\phi(x))\geq A_j(x)\, or\, B_{\psi(j)}(\phi(x))\geq A_j(x)\otimes l,$$
where $L\setminus 0=\{l\in L:l>0\}$.
\end{itemize}
\end{def1}
\begin{rem}\label{remark} (i) If $l=1,\mathcal{W}=\mathcal{P}_X\times \mathcal{P}_Y$ then $B_{\psi(j)}(\phi(x))\geq A_j(x)$ and the pair $(\phi,\psi)$ is an {\bf {FP}-map} from $(X,\mathcal{P}_X)$ to $(Y,\mathcal{P}_Y)$ given in \cite{mo}.\\\\
(ii) For a {{\text FP}-map} $(\phi,\psi,\mathcal{W}):(X,\mathcal{P}_X)\rightarrow (Y,\mathcal{P}_Y)$, where $(X,\mathcal{P}_X=\{A_{\xi(x)}:x\in X\})$ and $(Y,\mathcal{P}_Y=\{B_{\xi^{'}(y)}:y \in Y\})$, characterized by index functions $\xi:X\rightarrow J_1,\xi^{'}:Y\rightarrow J_2$, the diagram in Figure \ref{fig,3} commutes, i.e., $\xi^{'}\circ \phi=\psi\circ \xi$.
\begin{figure}
\begin{center}
\begin{tikzcd}[row sep=12ex, column sep=12ex]
X\arrow{r}{\xi}\arrow[swap]{d}{\phi}& J_1\arrow{d}{\psi}\\
Y\arrow{r}{\xi^{'}}&J_2
\end{tikzcd}
\end{center}
  \caption{Diagram for Remark \ref{remark}}
    \label{fig,3}
\end{figure}
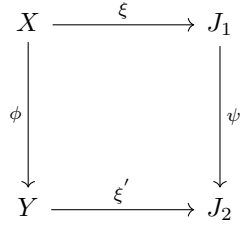
\end{rem}
\begin{pro} Spaces with $L$-fuzzy partitions alongwith their {\textit{FP}-maps} form a category.
\end{pro}
\textbf{Proof:} We only need to show that the composition of two \text{FP}-maps is again an {\text FP}-map, i.e, $(\phi_2,\psi_2,\mathcal{W}_2)\circ(\phi_1,\psi_1,\mathcal{W}_1)=(\phi_2\circ\phi_1,\psi_2\circ\psi_1,\mathcal{W}_2\cdot \mathcal{W}_1)$ is an {\text FP}-map. For which, let $(\phi_1,\psi_1,\mathcal{W}_1):(X,\mathcal{P}_X)\rightarrow (Y,\mathcal{P}_Y)$ and $(\phi_2,\psi_2,\mathcal{W}_2):(Y,\mathcal{P}_Y)\rightarrow (Z,\mathcal{P}_Z)$ be {\text FP}-maps, where $\mathcal{P}_X= \{A_j:j\in J_1\},\mathcal{P}_Y= \{B_{j^{'}}:j^{'}\in J_2\}$ and $\mathcal{P}_Z= \{C_{j^{''}}:j^{''}\in J_3\}$. Then there exist $l_1,l_2\in L\setminus 0$ such that  $ B_{\psi_1(j)}(\phi_1(x))\geq A_j(x)\otimes l_1$ and $C_{\psi_2(j^{'})}(\phi_2(y))\geq B_{j^{'}}(y)\otimes l_2,\, \forall\,x\in X,y\in Y,(A_j,B_{\psi_1(j)})\in \mathcal{W}_1,(B_{j^{'}},C_{\psi_2(j^{'})})\in \mathcal{W}_2$. Now, let $(A_j,C_{\psi_2(\psi_1(j))})\in \mathcal{W}_1\cdot \mathcal{W}_2$, $x\in X,j\in J_1$, where $\mathcal{W}_1\cdot \mathcal{W}_2$ is a classical composition of two relations. Then
\begin{eqnarray*}
C_{\psi_2\circ\psi_1(j)}(\phi_2\circ\phi_1(x))&=&C_{\psi_2(\psi_1(j))}(\phi_2(\phi_1(x)))\geq  B_{\psi_1(j)}(\phi_1(x))\otimes l_2\\
&\geq& (A_j(x)\otimes l_1)\otimes l_2=A_j(x)\otimes (l_1\otimes l_2).
\end{eqnarray*}
Thus the composition of two {\text FP}-maps is again an {\text FP}-map.\\\\
We shall denote by ${\bf FPS}$, the category of spaces with $L$-fuzzy partitions and their {{\text FP}-maps}.
\begin{pro} The category ${\bf FPS}$ has products.
\end{pro}
\textbf{Proof:} To show that $(X\times Y,\mathcal{P}_{X\times Y})$ is product of $(X,\mathcal{P}_X),(Y,\mathcal{P}_Y)\in|{\bf FPS}|$, we first need to show that $(p_1,h_1,\mathcal{W}_1)$ and $(p_2,h_2,\mathcal{W}_2)$ are ${\bf FPS}$-morphisms, where $\mathcal{P}_{X\times Y}=\{(A_j,B_{j,}):A_j\in \mathcal{P}_X,B_{j'}\in \mathcal{P}_Y,j\in J_1,j'\in J_2\} $ is an $L$-fuzzy partition of $X\times Y$. For which, let $p_1:X\times Y\rightarrow X,p_2:X\times Y\rightarrow Y$ and $h_1:J_1\times J_2\rightarrow J_1 ,h_2:J_1\times J_2\rightarrow J_2$ be projections associated with the cartesian products $X\times Y$ and $J_1\times J_2$, respectively, and $\mathcal{W}_1=\{(A_j\times B_{j^{'}},A_j):A_j\in \mathcal{P}_{X}, B_{j^{'}}\in \mathcal{P}_{ Y}\},\mathcal{W}_2=\{(A_j\times B_{j^{'}},B_{j^{'}}):A_j\in \mathcal{P}_{X}, B_{j^{'}}\in \mathcal{P}_{ Y}\}$. Then $A_j(x)\wedge B_{j^{'}}(y)\leq A_j(x)\Rightarrow (A_j\times B_{j^{'}})(x,y)=(A\times B)_{(j,j^{'})}(x,y)\otimes 1\leq A_j(x)=A_{h_1(j,j^{'})}(p_1(x,y))$ and similarly $(A_j\times B_{j^{'}})(x,y)\otimes 1\leq B_{h_2(j,j^{'})}(p_2(x,y))$, these projections are ${\bf FPS}$-morphisms. Now, let $(Z,\mathcal{P}_Z)\in |{\bf FPS}|$ and $(\phi_1,\psi_1,\mathcal{W}_1):(Z,\mathcal{P}_Z)\rightarrow (X,\mathcal{P}_X),(\phi_2,\psi_2,\mathcal{W}_2):(Z,\mathcal{P}_Z)\rightarrow (Y,\mathcal{P}_Y)$ be ${\bf FPS}$-morphisms. Further, we set $\phi:Z\rightarrow X\times Y,\psi:J_3\rightarrow J_1\times J_2$ such that $\phi(z)=(\phi_1(z),\phi_2(z)),\psi(j^{''})=(\psi_1(j^{''}),\psi_2(j^{''}))$ and $\mathcal{W}=\{(C_{j^{''}},A_{j}\times B_{j^{'}}):C_{j^{''}}\in \mathcal{P}_Z,A_{j}\times B_{j^{'}}\in \mathcal{P}_{X\times Y}\}$. Then for all $(C_{j^{''}},A_{j}\times B_{j^{'}})\in \mathcal{W}$, we have $(C_{j^{''}},A_{j})\in \mathcal{W}_1,(C_{j^{''}},B_{j^{'}}) \in \mathcal{W}_2$. Now, we obtain
\begin{eqnarray*}
C_{j^{''}}(z)\rightarrow (A\times B)_{\psi(j^{''})}(\phi(z))&=&C_{j^{''}}(z)\rightarrow (A\times B)_{(\psi_1(j^{''}),\psi_2(j^{''}))}(\phi_1(z),\phi_2(z))\\
&=&C_{j^{''}}(z)\rightarrow (A_{\psi_1(j^{''})}\times B_{\psi_2(j^{''})})(\phi_1(z),\phi_2(z))\\
&=&C_{j^{''}}(z)\rightarrow A_{\psi_1(j^{''})}(\phi_1(z))\wedge B_{\psi_2(j^{''})}(\phi_2(z))\\
&\geq& C_{j^{''}}(z)\rightarrow C_{j^{''}}(z)\otimes l_1\wedge C_{j^{''}}(z)\otimes l_2\\
&=&(C_{j^{''}}(z)\rightarrow C_{j^{''}}(z)\otimes l_1)\\
&&\wedge (C_{j^{''}}(z)\rightarrow C_{j^{''}}(z)\otimes l_2)\\
&\geq&l_1\otimes l_2.
\end{eqnarray*}
Thus $(\phi,\psi,\mathcal{W})$ is a ${\bf FPS}$-morphism and
the rest of the proof can be done easily.\\\\
Below, we recall the concept of $F^\uparrow$-transform from \cite{per,anan}.
\begin{def1}\label{DUFT}
Let $\mathcal{P}=\lbrace A_{j}:j\in J\rbrace$ be an $L$-fuzzy partition. Then the {\bf direct $F^\uparrow$-transform} of $f\in L^X$ is a collection of lattice elements $\{F^\uparrow_j[f]:j\in J\}$ and the $j$-th component of direct $F^\uparrow$-transform is given by
$$F_j^\uparrow[f]=\bigvee\limits_{x\in X}A_{j}(x)\otimes f(x),\, \forall\,j\in J,f\in L^X.$$
\end{def1}
\begin{pro}\label{3.1}
Let $\mathcal{P}=\lbrace A_{j}:j\in J\rbrace$ be an $L$-fuzzy partition. Then for all $a\in L,\textbf{a},f\in L^X$ and $\{f_i:i\in J\}\subseteq L^X$
\begin{itemize}
\item[(i)] $F_{j}^\uparrow[\textbf{a}]=a$;
\item[(ii)]  if $f\leq g$ then $F_{j}^\uparrow[f]\leq F_{j}^\uparrow[g]$;
\item[(iii)] $F_{j}^\uparrow[\textbf{a}\otimes f]=a\otimes F_{j}^\uparrow[f]$;
\item[(iv)] $F_{j}^\uparrow[\bigvee\limits_{i\in J}f_i]=\bigvee\limits_{i\in J}F_{j}^\uparrow[f_i]$; and
\item[(v)] $F_j^\uparrow[\bigwedge\limits_{i\in J}f_i]\leq\bigwedge\limits_{i\in J}F_{j}^\uparrow[f_i]$.
\end{itemize}
\end{pro}
\textbf{Proof:} Straightforward from Definitions \ref{DUFT}.\\\\
Now, we have the following.
\begin{pro}\label{1} $(\phi,\psi,\mathcal{W}):(X,\mathcal{P}_X)\rightarrow (Y,\mathcal{P}_Y)$ is a ${\bf FPS}$-morphism iff $ l \rightarrow F^\uparrow_{\psi(j)}[f]\geq F^\uparrow_{j}[ \overleftarrow{\phi}(f)]$ or $F^\uparrow_{\psi(j)}[f]\geq  F^\uparrow_{j}[ \overleftarrow{\phi}(f)]\otimes {l},\,\forall\,j\in J,l\in L\setminus 0,f\in L^Y,\mathcal{W}\subseteq \mathcal{P}_X\times \mathcal{P}_Y$.
\end{pro}
\textbf{Proof:} Let $(\phi,\psi,\mathcal{W}):(X,\mathcal{P}_X)\rightarrow (Y,\mathcal{P}_Y)$ be a ${\bf FPS}$-morphism. Then for all $f\in L^Y,j\in J,\mathcal{W}\subseteq \mathcal{P}_X\times \mathcal{P}_Y$,
\begin{eqnarray*}
F^\uparrow_{j}[ \overleftarrow{\phi}(f)]&=&\bigvee\limits_{x\in X}A_{j}(x)\otimes \overleftarrow{\phi}(f)(x)\\
&=&\bigvee\limits_{x\in X}A_{j}(x)\otimes f(\phi(x))\\
&\leq&\bigvee\limits_{x\in X}(l\rightarrow B_{\psi(j)}(\phi(x)))\otimes f(\phi(x))\\
&\leq&\bigvee\limits_{y\in Y}(l\rightarrow B_{\psi(j)}(y))\otimes f(y)\\
&\leq&\bigvee\limits_{y\in Y}(l\rightarrow B_{\psi(j)}(y)\otimes f(y))\\
&\leq& l\rightarrow \bigvee\limits_{y\in Y}B_{\psi(j)}(y)\otimes f(y)\\
&=& {l}\rightarrow F^\uparrow_{\psi(j)}[f]. 
\end{eqnarray*}
Thus $ l \rightarrow F^\uparrow_{\psi(j)}[f]\geq F^\uparrow_{j}[ \overleftarrow{\phi}(f)]$ or $F^\uparrow_{\psi(j)}[f]\geq  F^\uparrow_{j}[ \overleftarrow{\phi}(f)]\otimes {l}$.\\\\
Conversely, let $\phi:X\rightarrow Y:\psi,J_1\rightarrow J_2$ and $\mathcal{W}\subseteq \mathcal{P}_X\times \mathcal{P}_Y, \textit{\textbf{dom}}(\mathcal{W})=\mathcal{P}_X$ such that $l \rightarrow F^\uparrow_{\psi(j)}[f]\geq F^\uparrow_{j}[ \overleftarrow{\phi}(f)]$ or $F^\uparrow_{\psi(j)}[f]\geq  F^\uparrow_{j}[ \overleftarrow{\phi}(f)]\otimes {l},\,\forall\,j\in J,f\in L^Y$. Then for all $x\in X,j\in J$ and $(A_j,B_{\psi(j)})\in \mathcal{W}$,
\begin{eqnarray*}
B_{\psi(j)}(\phi(x))&=&\bigvee\limits_{y\in Y}B_{\psi(j)}(y)\otimes1_{\{\phi(x)\}}(y)\\
&=&F^\uparrow_{\psi(j)}[1_{\{\phi(x)\}}]\\
&\geq&  F^\uparrow_{j}[\overleftarrow{\phi}(1_{\{\phi(x)\}})]\otimes {l}\\
&=&\bigvee\limits_{z\in X}(A_{j}(z)\otimes \overleftarrow{\phi}(1_{\{\phi(x)\}})(z))\otimes l\\
&=&\bigvee\limits_{z\in X}(A_{j}(z)\otimes 1_{\{\phi(x)\}}(\phi(z)))\otimes l\\
&=&A_{j}(x)\otimes l.
\end{eqnarray*}
Thus $(\phi,\psi,\mathcal{W}):(X,\mathcal{P}_X)\rightarrow (Y,\mathcal{P}_Y)$ is a ${\bf FPS}$-morphism.
\begin{pro}\label{11} If $(\phi,\psi,\mathcal{W}):(X,\mathcal{P}_X)\rightarrow (Y,\mathcal{P}_Y)$ is a ${\bf FPS}$-morphism. Then $ l \rightarrow F^\uparrow_{\psi(j)}[ \overrightarrow{\phi}(f)]\geq F^\uparrow_{j}[f]$ or $F^\uparrow_{\psi(j)}[ \overrightarrow{\phi}(f)]\geq F^\uparrow_{j}[f]\otimes {l},\,\forall\,j\in J,l\in L\setminus 0,f\in L^X,\mathcal{W}\subseteq \mathcal{P}_X\times \mathcal{P}_Y$.
\end{pro}
\textbf{Proof:} Similar to that of Proposition \ref{1}.\\\\
In the following, we introduce the concept of an $L$-fuzzy approximation space from categorical point of view.
\begin{def1}
Let $(X,R_X)$ and $(X,R_Y)$ be two $L$-fuzzy approximation space. Then $\phi:(X,R_X)\rightarrow(X,R_Y)$ is an {\bf order preserving function}, if
\begin{itemize}
\item[(i)] $\phi:X\rightarrow Y$ is a function; and
\item[(ii)] $\exists\,l\in L\setminus 0$ such that $l\rightarrow R_Y(\phi(x),\phi(y))\geq R_X(x,y)$ or $R_Y(\phi(x),\phi(y))\geq R_X(x,y)\otimes l,\,\forall\,x,y\in X$.
\end{itemize}
\end{def1}
\begin{rem} If $l=1$. Then $R_Y(\phi(x),\phi(y))\geq R_X(x,y)$ and $\phi$ is an order preserving function between two $L$-fuzzy approximation spaces given in \cite{lai}.
\end{rem}
\begin{pro}\label{4.3} $L$-fuzzy approximation spaces alongwith their order preserving functions form a category.
\end{pro}
\textbf{Proof:} We only need to show that the composition of two order preserving functions is again an order preserving function. For which, let $\phi_1:(X,R_X)\rightarrow (Y,R_Y)$ and $\phi_2:(Y,R_Y)\rightarrow (Z,R_Z)$ be order preserving functions, i.e., $\phi_1:X\rightarrow Y,\phi_2:Y\rightarrow Z$ are functions and there exist $l_1,l_2\in L\setminus 0$ such that  $R_Y(\phi_1(x),\phi_1(z))\geq R_{X}(x,z)\otimes {l_1},R_Z(\phi_2(y),\phi_2(u))\geq R_{Y}(y,u)\otimes {l_2},\,\forall\,x,z\in X,y,u\in Y$. Now, let $x,z\in X$. Then 
\begin{eqnarray*}
R_Z((\phi_2\circ\phi_1)(x),(\phi_2\circ\phi_1)(z))&=&R_Z(\phi_2(\phi_1(x)),\phi_2(\phi_1(z)))\\
&\geq& R_Y(\phi_1(x),\phi_1(z))\otimes {l_2}\\
&\geq &(R_{X}(x,z)\otimes {l_1})\otimes l_2\\
&=& R_{X}(x,z)\otimes ( l_1\otimes l_2).
\end{eqnarray*}
Thus $\phi_2\circ\phi_1:(X,\Upsilon_X)\rightarrow (Z,\Upsilon_Z)$ is an order preserving function.\\\\
We shall denote by ${\bf FAS}$, the category of $L$-fuzzy approximation spaces and their order preserving functions.
\begin{pro}\label{ua} $\phi:(X,{R}_X)\rightarrow (Y,R_Y)$ is a ${\bf FAS}$-morphism iff ${{l}}\rightarrow \overleftarrow{\phi}(\overline{R}_Y(f))(x)\geq \overline{R}_X(\overleftarrow{\phi}(f))(x)\,or\,\overleftarrow{\phi}(\overline{R}_Y(f))(x)\geq  \overline{R}_X(\overleftarrow{\phi}(f))(x)\otimes {{l}},\,\forall\,x\in X,f\in L^Y,{l}\in L\setminus {0}$.
\end{pro}
\textbf{Proof:} Similar to that of Proposition \ref{1}.
\begin{rem}\label{op}
If $(\phi,\psi,\mathcal{W}):(X,\mathcal{P}_X)\rightarrow (Y,\mathcal{P}_Y)$ is a ${\bf FPS}$-morphism. Then $\phi:(X,R_{\mathcal{P}_X})\rightarrow (Y,R_{\mathcal{P}_Y})$ is a ${\bf FAS}$-morphism.
\end{rem}
\begin{pro}\label{F_1}
Let $F_1: {\bf FPS}\rightarrow {\bf FAS}$ be a function such that for all $(X,\mathcal{P}_X)\in|{\bf FPS}|,\,F_1(X,\mathcal{P}_X)=(X,R_{\mathcal{P}_X})$ and for every ${\bf FPS}$-morphism $(\phi,\psi,\mathcal{W}):(X,\mathcal{P}_X)\rightarrow (Y,\mathcal{P}_Y)$, $F_1(\phi,\psi,\mathcal{W}):(X,R_{\mathcal{P}_X})\rightarrow (Y,R_{\mathcal{P}_Y})$ such that $F_1(\phi,\psi,\mathcal{W})=\phi$. Then $F_1$ is a functor.
\end{pro} 
\textbf{Proof:} Follows from Proposition \ref{op}.
\section{Categorical view of $L$-fuzzy partitions and $L$-fuzzy closure systems}
In this section, we recall the concept of $L$-fuzzy closure system space from \cite{fan}. Further, we show that a relationship between the category of spaces with $L$-fuzzy partitions and the category $L$-fuzzy closure system spaces. Now, we begin with the following.
\begin{def1} \label{FCS}
A function $\Upsilon : L^X \rightarrow L$ is called an $L$-fuzzy closure system on $X$ if it satisfies the following conditions:
\begin{itemize}
\item[(i)] $\Upsilon(1_X)=1$; and
\item[(ii)] $\Upsilon(\bigwedge\limits_{i\in J}f_i)\geq \bigwedge\limits_{i\in J}\Upsilon(f_i),\,\forall\,\{f_i,i\in J\}\subseteq L^X$.
\end{itemize}
The pair $(X,\Upsilon)$ is called an {\bf $L$-fuzzy closure system space}. 
\end{def1}
An $L$-fuzzy closure system space $(X,\Upsilon)$ is called 
\begin{itemize}
\item {\bf enriched}, if $\Upsilon(\textbf{a}\rightarrow f)\geq \Upsilon(f) ,\,\forall\,\textbf{a},f\in L^X$; and
\item {\bf strong}, if $\Upsilon(\textbf{a} \otimes f)\geq \Upsilon(f) ,\,\forall\,\textbf{a},f\in L^X$.
\end{itemize}
Next, we introduce the concept of an $L$-fuzzy closure system space from categorical point of view by using slightly different continuous function. Now, we begin with the following.
\begin{def1} \label{cc}
Let $(X,\Upsilon_X)$ and $(Y,\Upsilon_Y)$ be $L$-fuzzy closure system spaces. Then $\phi:(X,\Upsilon_X)\rightarrow (Y,\Upsilon_Y)$ is a {\bf continuous function} if
\begin{itemize}
    \item[(i)] $\phi:X\rightarrow Y$ is a function; and \item[(ii)] $\exists\, l\in L\setminus 0$ such that $l\rightarrow \Upsilon_X(\overleftarrow{\phi}(f))\geq \Upsilon_{Y}(f)\,or\,\Upsilon_X(\overleftarrow{\phi}(f))\geq \Upsilon_{Y}(f)\otimes {l},\,\forall\,f\in L^Y$.
 \end{itemize}
\end{def1}
\begin{rem} If $l=1$ then $\Upsilon_X(\overleftarrow{\phi}(f))\geq \Upsilon_{Y}(f)$ and $\phi$ is a continuous function between $L$-fuzzy closure system spaces given in \cite{fan}.
\end{rem}
\begin{pro}\label{clc} $L$-fuzzy closure system spaces alongwith their continuous functions form a category.
\end{pro}
\textbf{Proof:} We only need to show that the composition of two continuous functions is again a continuous function. For which, let $\phi_1:(X,\Upsilon_X)\rightarrow (Y,\Upsilon_Y)$ and $\phi_2:(Y,\Upsilon_Y)\rightarrow (Z,\Upsilon_Z)$ be continuous functions, i.e., $\phi_1:X\rightarrow Y,\phi_2:Y\rightarrow Z$ are functions and there exist $l_1,l_2\in L\setminus 0$ such that  $\Upsilon_X(\overleftarrow{\phi_1}(f))\geq \Upsilon_{Y}(f)\otimes{l_1},\Upsilon_Y(\overleftarrow{\phi_2}(g))\geq \Upsilon_{Z}(g)\otimes {l_2},\,\forall\,f\in L^Y,\,g\in L^Z$. Now, let $g\in L^Z$. Then 
\begin{eqnarray*}
\Upsilon_X(\overleftarrow{\phi_2\circ \phi_1}(g))&=&\Upsilon_X((\overleftarrow{\phi_1}\circ \overleftarrow{\phi_2})(g))\\
&=&\Upsilon_X(\overleftarrow{\phi_1}( \overleftarrow{\phi_2}(g)))\\
&\geq &\Upsilon_Y( \overleftarrow{\phi_2}(g))\otimes l_1\\
&\geq& (\Upsilon_Z(g)\otimes l_2)\otimes l_1\\
&=&\Upsilon_Z(g)\otimes ( l_2\otimes l_1).
\end{eqnarray*}
Thus $\phi_2\circ\phi_1:(X,\Upsilon_X)\rightarrow (Z,\Upsilon_Z)$ is a continuous function.\\\\
We shall denote by ${\bf FCSS}$, the category of $L$-fuzzy closure system spaces and their continuous functions. Further, we shall denote by ${\bf EFCSS}$, the full subcategory of ${\bf FCSS}$ with objects as enriched $L$-fuzzy closure system spaces and their continuous functions.
\begin{pro} The $L$-fuzzy closure system space $(\emptyset,\Upsilon_{\emptyset})\in|{\bf FCSS}|$ is an initial object.
\end{pro}
\textbf{Proof:} Let $(\emptyset,\Upsilon_{\emptyset}),(X,\Upsilon_X) \in |{\bf FCSS}|$ and $\phi:\emptyset\rightarrow X$ be a function. Then $\phi$ is a unique function (as uniqueness condition for empty function vacuously true). Therefore the condition (ii) in Definition \ref{cc} hold for $\phi$. Thus $(\emptyset,\Upsilon_{\emptyset})\in |{\bf FCSS}|$ is an initial object.
\begin{pro} The $L$-fuzzy closure system space $(\{*\},\Upsilon_{\{*\}})\in|{\bf FCSS}|$ is a final object, where $\{*\}$ is a singleton. 
\end{pro}
\textbf{Proof:} Let $(\{*\},\Upsilon_{\{*\}}),(X,\Upsilon_X)\in|{\bf FCSS}|$ and $L^{{\{*\}}}=\{0_{{\{*\}}},1_{{\{*\}}}\}$. Next, we have to show that there exist a unique ${\bf FCSS}$-morphism from $(X,\Upsilon_X)$ to $(\{*\},\Upsilon_{\{*\}})$. Further, let $\phi:X\rightarrow {\{*\}}$ is a function, which is a unique function. Clearly, the condition (ii) in Definition \ref{cc} hold for $\phi$. Thus $({\{*\}},\Upsilon_{{\{*\}}})\in|{\bf FCSS}|$ is a final object.\\\\
In the following proposition, we show that ${\bf FCSS}$ is a semiadditive category. For this, we take some assumptions. Let $X$ be a commutative monoid with respect to addition (`0' is the additive identity) and $f\in L^X$ is defined as $f(x+y)=f(x)\wedge f(y)$. Further, let the collection $\mathcal{A}=\{\phi|\phi:X\rightarrow Y\,\text{is\, a\, function} \}$ be equipped with the structure of a commutative monoid with operation `+', where zero function is an additive identity. For any $\phi_1,\phi_2\in \mathcal{A},\phi_1+\phi_2\in \mathcal{A}$ and we define $(\phi_1+\phi_2)(x)=\phi_1(x)+\phi_2(x)$. Next, according to Zadeh's extenstion principle the collection can be extended to the collection of operators $\mathcal{B}=\{\overleftarrow{\phi}:L^Y\rightarrow L^X|\overleftarrow{\phi}(f)=f\circ \phi,\,\forall\,f\in L^Y\}$. Next, for $\phi_1+\phi_2:X\rightarrow Y$, we have $\overleftarrow{\phi_1+\phi_2}:L^Y\rightarrow L^X$ such that 
$\overleftarrow{\phi_1+\phi_2}(f)(x)=(f\circ(\phi_1+\phi_2))(x),\,\forall\,x\in X,f\in L^Y$.\\
As we know that for all $x\in X,f\in L^Y,\,\overleftarrow{\phi_1+\phi_2}(f)(x)=(f\circ(\phi_1+\phi_2))(x)=f((\phi_1+\phi_2)(x))=f(\phi_1(x)+\phi_2(x))=f(\phi_1(x))\wedge f(\phi_2(x))=(f\circ \phi_1)(x)\wedge (f\circ \phi_2)(x)=((f\circ \phi_1)\wedge (f\circ \phi_2))(x)=(\overleftarrow{\phi_1}(f)\wedge \overleftarrow{\phi_2}(f))(x).$ Thus $\overleftarrow{\phi_1+\phi_2}(f)=\overleftarrow{\phi_1}(f)\wedge \overleftarrow{\phi_2}(f)$.
\begin{pro}
The semiadditive structure on homset in ${\bf FCSS}$ is given by $\phi_1+\phi_2$. Here zero function serves as additive identity. 
\end{pro}
\textbf{Proof:} Let $\phi_1,\phi_2:(X,\Upsilon_X)\rightarrow (Y,\Upsilon_Y)$ be ${\bf FCSS}$-morphisms. To show $\phi_1+\phi_2:(X,\Upsilon_X)\rightarrow (Y,\Upsilon_Y)$ be a ${\bf FCSS}$-morphism, let $(\phi_1+\phi_2)(x)=\phi_1(x)+\phi_2(x)$ and $f(x+y)=f(x)\wedge f(y),\,\forall \,x\in X,f\in L^Y$. Now, $\Upsilon_X(\overleftarrow{\phi_1+\phi_2}(f))=\Upsilon_X(\overleftarrow{\phi_1}(f)\wedge \overleftarrow{\phi_2}(f))\geq \Upsilon_X(\overleftarrow{\phi_1}(f))\wedge \Upsilon_X(\overleftarrow{\phi_2}(f))\geq (\Upsilon_Y(f)\otimes l)\wedge (\Upsilon_Y(f)\otimes l)=\Upsilon_Y(f)\otimes l$. Thus $\phi_1+\phi_2:(X,\Upsilon_X)\rightarrow (Y,\Upsilon_Y)$ is a ${\bf FCSS}$-morphism with the zero function as additive identity. Also, composition distributes over `+', i.e., for $\phi_1:X\rightarrow Y,\phi_2,\phi_3:Y\rightarrow Z,x\in X,$ $((\phi_2+\phi_3)\circ \phi_1)(x)= (\phi_2+\phi_3)(\phi_1(x))=\phi_2(\phi_1(x))+\phi_3(\phi_1(x))=(\phi_2\circ\phi_1)(x)+(\phi_3\circ\phi_1)(x)=(\phi_2\circ\phi_1+\phi_3\circ\phi_1)(x).$ Thus $(\phi_2+\phi_3)\circ \phi_1=\phi_2\circ\phi_1+\phi_3\circ\phi_1$, i.e., composition distributes over `+'.
\\\\
In the following, we demonstrate that an upper $F$-transform determines an $L$-fuzzy closure system uniquely.
\begin{pro} \label{3.3}
Let $(X,\mathcal{P})$ be a space with $L$-fuzzy partition $\mathcal{P}=\lbrace A_{\xi(x)}:x\in X\rbrace$ and
$$ \Upsilon_{\mathcal{P}}(f)=\bigwedge\limits_{x\in X}(F_{\xi(x)}^\uparrow[f]\rightarrow f(x)),\,\forall\, f\in L^X.$$
where $\xi:X\rightarrow J$ is a function such that $\xi(x)$ is the unique element of $J$ with $x\in core(A_{\xi(x)})$.
 Then $(X,\Upsilon_{\mathcal{P}})$ is enriched and strong $L$-fuzzy closure system space.
\end{pro}
\textbf{Proof:} (i) From Propositions \ref{3.1} and \ref{3.3}
$$\Upsilon_{\mathcal{P}}(1_X)=\bigwedge\limits_{x\in X}(F_{\xi(x)}^\uparrow[1_X]\rightarrow 1_X(x))=\bigwedge\limits_{x\in X}(1\rightarrow 1)=1.$$
(ii) Let $\{f_i,i\in J\}\subseteq L^X$. Then from Propositions \ref{3.1} and \ref{3.3}
\begin{eqnarray*}
\Upsilon_{\mathcal{P}}(\bigwedge\limits_{i\in J}f_i)&=&\bigwedge\limits_{x\in X}(F_{\xi(x)}^\uparrow[\bigwedge\limits_{i\in J}f_i]\rightarrow \bigwedge\limits_{i\in J}f_i(x))\\
&\geq&\bigwedge\limits_{x\in X}\bigwedge\limits_{i\in J}(\bigwedge\limits_{i\in J}F_{\xi(x)}^\uparrow[f_i]\rightarrow f_i(x))\\
&\geq&\bigwedge\limits_{x\in X}\bigwedge\limits_{i\in J}(F_{\xi(x)}^\uparrow[f_i]\rightarrow f_i(x))\\
&=&\bigwedge\limits_{i\in J}\Upsilon_{\mathcal{P}}(f_i).
\end{eqnarray*}
(iii) Let $\textbf{a},f\in L^X$. Then from Propositions \ref{3.1} and \ref{3.3} 
\begin{eqnarray*}
\Upsilon_{\mathcal{P}}(\textbf{a}\rightarrow f)&=&\bigwedge\limits_{x\in X}(F_{\xi(x)}^\uparrow[\textbf{a}\rightarrow f]\rightarrow (\textbf{a}\rightarrow f)(x))\\
&=&\bigwedge\limits_{x\in X}(F_{\xi(x)}^\uparrow[\textbf{a}\rightarrow f]\rightarrow ({a}\rightarrow f(x)))\\
&=&\bigwedge\limits_{x\in X}(a\otimes F_{\xi(x)}^\uparrow[\textbf{a}\rightarrow f]\rightarrow f(x))\\
&=&\bigwedge\limits_{x\in X} (F_{\xi(x)}^\uparrow[\textbf{a}\otimes(\textbf{a}\rightarrow f)]\rightarrow f(x))\\
&\geq&\bigwedge\limits_{x\in X}( F_{\xi(x)}^\uparrow[f]\rightarrow f(x))\\
&=&\Upsilon_{\mathcal{P}}(f).
\end{eqnarray*}
(iv) Let $\textbf{a},f\in L^{X}$. Then from Propositions \ref{3.1} and \ref{3.3} 
\begin{eqnarray*}
\Upsilon_{\mathcal{P}}(\textbf{a}\otimes f)&=&\bigwedge\limits_{x\in X}(F_{\xi(x)}^\uparrow[\textbf{a}\otimes f]\rightarrow (\textbf{a}\otimes f)(x))\\
&=&\bigwedge\limits_{x\in X}(F_{\xi(x)}^\uparrow[\textbf{a}\otimes f]\rightarrow ({a}\otimes f(x)))\\
&=&\bigwedge\limits_{x\in X}(a\otimes F_{\xi(x)}^\uparrow[f]\rightarrow a\otimes f(x))\\
&\geq& \bigwedge\limits_{x\in X} (F_{\xi(x)}^\uparrow[f]\rightarrow f(x))\\
&=&\Upsilon_{\mathcal{P}}(f).
\end{eqnarray*}
\begin{pro}\label{2}
If $(\phi,\psi,\mathcal{W}):(X,\mathcal{P}_X)\rightarrow (Y,\mathcal{P}_Y)$ is a ${\bf FPS}$-morphism. Then $\phi:(X,\Upsilon_{\mathcal{P}_X})\rightarrow (Y,\Upsilon_{\mathcal{P}_Y})$ is a ${\bf FCSS}$-morphism.
\end{pro}
\textbf{Proof:} Let $f\in L^Y$ and  $(\phi,\psi,\mathcal{W}):(X,\mathcal{P}_X)\rightarrow (Y,\mathcal{P}_Y)$ be a ${\bf FPS}$-morphism. Then from Proposition \ref{3.3}
\begin{eqnarray*}
\Upsilon_{\mathcal{P}_X}(\overleftarrow{\phi}(f))&=&\bigwedge\limits_{x\in X}(F_{\xi(x)}^\uparrow[\overleftarrow{\phi}(f)]\rightarrow \overleftarrow{\phi}(f)(x))\\
&\geq& \bigwedge\limits_{x\in X}((l\rightarrow F_{\psi(\xi(x))}^\uparrow[f])\rightarrow f(\phi(x)))\\
&\geq&\bigwedge\limits_{x\in X}l\otimes (F_{\psi(\xi(x))}^\uparrow[f]\rightarrow f(\phi(x)))\\
&\geq& l\otimes \bigwedge\limits_{x\in X}(F_{\psi(\xi(x))}^\uparrow[f]\rightarrow f(\phi(x)))\\
&=&l\otimes \bigwedge\limits_{x\in X}(F_{\xi^{'}(\phi(x))}^\uparrow[f]\rightarrow f(\phi(x)))\\
&\geq& l\otimes \bigwedge\limits_{y\in Y}(F_{\xi^{'}(y)}^\uparrow[f]\rightarrow f(y))\\
&=&l\otimes \Upsilon_{\mathcal{P}_Y}(f).
\end{eqnarray*}
Thus $\phi:(X,\Upsilon_{\mathcal{P}_X})\rightarrow (Y,\Upsilon_{\mathcal{P}_Y})$ is a ${\bf FCSS}$-morphism.
\begin{pro}\label{F_3}
Let $F_3:{\bf FPS}\rightarrow {\bf FCSS}$ be a function such that for all $(X,\mathcal{P}_X)\in|{\bf FPS}|,\,F_3(X,\mathcal{P}_X)=(X,\Upsilon_{\mathcal{P}_X})$ and for every ${\bf FPS}$-morphism $(\phi,\psi,\mathcal{W}):(X,\mathcal{P}_X)\rightarrow (Y,\mathcal{P}_Y)$, $F_3(\phi,\psi,\mathcal{W}):(X,\Upsilon_{\mathcal{P}_X})\rightarrow (Y,\Upsilon_{\mathcal{P}_Y})$ such that $F_3(\phi,\psi,\mathcal{W})=\phi$. Then $F_3$ is a functor.
\end{pro} 
\textbf{Proof:} Follows from Propositions \ref{3.3} and \ref{2}.\\\\
Next, we introduce the category of $L$-fuzzy closure spaces. We begin with the following.
\begin{def1} 
Let $(X,c_X)$ and $(Y,c_Y)$ be $L$-fuzzy closure spaces. Then $\phi:(X,c_X)\rightarrow (Y,c_Y)$ is a {\bf continuous function}, if
\begin{itemize}
\item[(i)] $\phi:X\rightarrow Y$ is function; and
\item[(ii)] $\exists\,{l}\in L\setminus {0}$ such that ${l}\rightarrow\overleftarrow{\phi}(c_{Y}(f))(x)\geq c_X(\overleftarrow{\phi}(f))(x)\,or\,\overleftarrow{\phi}(c_{Y}(f))(x)\geq c_X(\overleftarrow{\phi}(f))(x)\otimes {{l}},\,\forall\,x\in X,f\in L^Y$.
\end{itemize}
\end{def1}
\begin{rem}
If ${\textit{l}}=1$. Then $\overleftarrow{\phi}(c_{Y}(f))\geq c_X(\overleftarrow{\phi}(f))$ and $\phi$ is a continuous function between two $L$-fuzzy closure spaces given in \cite{bel}.
\end{rem}
\begin{pro}
$L$-fuzzy closure spaces alongwith their continuous functions form a category.
\end{pro} 
\textbf{Proof:} Similar to that of Proposition \ref{clc}.\\\\
We shall denote by ${\bf FCS}$, the category of $L$-fuzzy closure spaces and their continuous functions. Further, we shall denote by ${\bf SFCS}$, the full subcategory of ${\bf FCS}$ with objects as strong $L$-fuzzy closure  spaces and their continuous functions.\\\\
In the following proposition, we present a functorial relationship among the above-introduced categories.
\begin{pro}\label{f}
There exist functors such that the diagram in Figure \ref{fig,2} commutes.
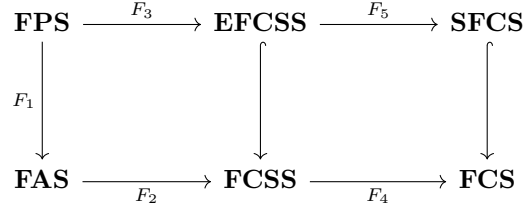
\begin{figure}
\begin{center}
\begin{tikzcd}[row sep=10ex, column sep=10ex] 
{\bf FPS}{{\arrow{r}{F_3}}}\arrow[swap]{d}{F_1}& {\bf EFCSS}{\arrow{r}{F_5}}\arrow[d, hook]& {\bf SFCS}\arrow[d, hook]\\
{\bf FAS}\arrow[r,"F_{2}"']& {\bf FCSS}\arrow[r,"F_{4}"']&{\bf FCS}
\end{tikzcd}
\end{center}
    \caption{Diagram for Proposition \ref{f}}
    \label{fig,2}
\end{figure}
\end{pro}
\textbf{Proof:} (1) From Proposition \ref{F_1}, $F_1:{\bf FPS}\rightarrow {\bf FAS}$ is a functor.\\\\
(2) For a ${\bf FAS}$-morphism  $\phi:(X,R_X)\rightarrow (Y,R_Y)$, we define a function $F_2: {\bf FAS}\rightarrow {\bf FCSS}$ such that
\[F_2(X,R_{X})=(X,\Upsilon_{R_X}),\,F_2(Y,R_Y)=(Y,\Upsilon_{R_Y}),\,F_2(\phi)=\phi\]
\[\Upsilon_{R_X}(f)=\bigwedge\limits_{x\in X}(\overline{R}_X(f)(x)\rightarrow f(x)),\,\forall\,f\in L^X,\]
\[\Upsilon_{R_Y}(g)=\bigwedge\limits_{y\in Y}(\overline{R}_Y(g)(y)\rightarrow g(y)),\,\forall\,g\in L^Y.\]
Then for all $g\in L^Y$
\begin{eqnarray*}
\Upsilon_{R_X}(\overleftarrow{\phi}( g))&=&\bigwedge\limits_{x\in X}(\overline{R}_X(\overleftarrow{\phi}(g))(x)\rightarrow \overleftarrow{\phi}(g)(x))\\
&=&\bigwedge\limits_{x\in X}(\overline{R}_X(\overleftarrow{\phi}(g))(x)\rightarrow g(\phi(x)))\\
&\geq&\bigwedge\limits_{x\in X}((l\rightarrow \overleftarrow{\phi}(\overline{R}_Y(g))(x))\rightarrow g(\phi(x)))\\
&\geq& \bigwedge\limits_{x\in X}l\otimes (\overleftarrow{\phi}(\overline{R}_Y(g))(x)\rightarrow g(\phi(x)))\\
&\geq&  l\otimes \bigwedge\limits_{x\in X}(\overline{R}_Y(g)(\phi(x))\rightarrow g(\phi(x)))\\
&\geq& l\otimes \bigwedge\limits_{y\in Y}(\overline{R}_Y(g)(y)\rightarrow g(y))\\
&=& l\otimes \Upsilon_{R_Y}(g).
\end{eqnarray*}
Thus $\phi:(X,\Upsilon_{R_X})\rightarrow (Y,\Upsilon_{R_Y})$ is a ${\bf FCSS}$-morphism and $F_2$ is a functor. For a ${\bf FCSS}$-morphism $\phi:(X,\Upsilon_X)\rightarrow (Y,\Upsilon_Y)$, we define a function $F^{-1}_2:{\bf FCSS}\rightarrow {\bf FAS}$ such that
\[F^{-1}_2(X,\Upsilon_{X})=(X,R_{\Upsilon_{X}}),\,\,F^{-1}_2(Y,\Upsilon_Y)=(Y,R_{\Upsilon_{Y}}),F^{-1}_2(\phi)=\phi.\]
\[R_{\Upsilon_{X}}(x,z)=\bigwedge\limits_{f\in L^X}(\Upsilon_X(f)\rightarrow (f(x)\rightarrow f(z))),\,\forall\,x,z\in X\]
\[R_{\Upsilon_{Y}}(y,u)=\bigwedge\limits_{g\in L^Y}(\Upsilon_Y(g)\rightarrow (g(y)\rightarrow g(u))),\,\forall\,y,u\in Y.\]
Then for all $x,z\in X$,
\begin{eqnarray*}
R_{\Upsilon_X}(x,z)&=&\bigwedge\limits_{f\in L^X}(\Upsilon_X(f)\rightarrow(f(x)\rightarrow f(z)))\\
&\leq& \bigwedge\limits_{g\in L^Y}(\Upsilon_X(\overleftarrow{\phi}(g))\rightarrow(\overleftarrow{\phi}(g)(x)\rightarrow \overleftarrow{\phi}(g)(z)))\\
&\leq& \bigwedge\limits_{g\in L^Y}(l\otimes \Upsilon_Y(g)\rightarrow(\overleftarrow{\phi}(g)(x)\rightarrow \overleftarrow{\phi}(g)(z)))\\
&=&\bigwedge\limits_{g\in L^Y}(l\rightarrow (\Upsilon_Y(g)\rightarrow(g(\phi(x))\rightarrow g(\phi(z)))))\\
&=&l\rightarrow \bigwedge\limits_{g\in L^Y}(\Upsilon_Y(g)\rightarrow(g(\phi(x))\rightarrow g(\phi(z))))\\
&=& l\rightarrow R_{\Upsilon_Y}(\phi(x),\phi(z)).
\end{eqnarray*}
Thus $\phi:(X,R_{\Upsilon_X})\rightarrow (Y,R_{\Upsilon_Y})$ is a ${\bf FAS}$-morphism and $F^{-1}_2$ is a functor and it is clear that $F_2$ and $F^{-1}_2$ are inverse functors.\\\\
(3) From proposition \ref{F_3}, $F_3:{\bf FPS}\rightarrow {\bf EFCSS}$ is a functor.\\\\
(4) For a ${\bf FCSS}$-morphism $\phi:(X,\Upsilon_X)\rightarrow (Y,\Upsilon_Y)$, we define a function $F_4:{\bf FCSS}\rightarrow {\bf FCS}$ such that
\[F_4(X,\Upsilon_{X})=(X,c_{\Upsilon_{X}}),\,F_4(Y,\Upsilon_Y)=(Y,c_{\Upsilon_{Y}}),\,F_4(\phi)=\phi\]
\[c_{\Upsilon_{X}}(f)(x)=\bigwedge\limits_{g\in L^X}(\Upsilon_X(g)\otimes \bigwedge\limits_{z\in X}(f(z)\rightarrow g(z))\rightarrow g(x)),\,\forall\,x\in X,f\in L^X\]
\[c_{\Upsilon_{Y}}(h)(y)=\bigwedge\limits_{k\in L^Y}(\Upsilon_Y(k)\otimes \bigwedge\limits_{u\in Y}(h(u)\rightarrow k(u))\rightarrow k(y)),\,\forall\,y\in Y,h\in L^Y.\]
Now, for all $x\in X$ and $h\in L^Y$,
\begin{eqnarray*}
c_{\Upsilon_{X}}(\overleftarrow{\phi}(h))(x)&=&\bigwedge\limits_{f\in L^X}(\Upsilon_X(f)\otimes \bigwedge\limits_{z\in Z}(\overleftarrow{\phi}(h)(z)\rightarrow f(z))\rightarrow f(x))\\
&\leq&\bigwedge\limits_{k\in L^Y}(\Upsilon_X(\overleftarrow{\phi}(k))\otimes \bigwedge\limits_{z\in Z}(\overleftarrow{\phi}(h)(z)\rightarrow \overleftarrow{\phi}(k)(z))\rightarrow \overleftarrow{\phi}(k)(x))\\
&\leq&\bigwedge\limits_{k\in L^Y}((\Upsilon_Y(k)\otimes l)\otimes \bigwedge\limits_{z\in Z}(h(\phi(z))\rightarrow k(\phi(z)))\rightarrow k(\phi(x)))\\
&\leq&\bigwedge\limits_{k\in L^Y}((\Upsilon_Y(k)\otimes l)\otimes \bigwedge\limits_{y\in Y}(h(y)\rightarrow k(y))\rightarrow k(\phi(x)))\\
&=&l\rightarrow\bigwedge\limits_{k\in L^Y}(\Upsilon_Y(k)\otimes \bigwedge\limits_{y\in Y}(h(y)\rightarrow k(y))\rightarrow k(\phi(x)))\\
&=&l\rightarrow c_Y(h)(\phi(x))\\
&=&l\rightarrow\overleftarrow{\phi}( c_Y(h))(x).
\end{eqnarray*}
Thus $\phi:(X,c_{\Upsilon_{X}})\rightarrow (Y,c_{\Upsilon_{Y}})$ is a ${\bf FCS}$-morphism and $F_4$ is a functor. Also, for a ${\bf FCS}$-morphism $\phi:(X,c_X)\rightarrow (Y,c_Y)$, we define a function $F^{-1}_4: {\bf FCS}\rightarrow {\bf FCSS}$ such that
\[F^{-1}_4(X,c_{X})=(X,\Upsilon_{c_X}),\,F^{-1}_4(Y,c_Y)=(Y,\Upsilon_{c_Y}),\,F^{-1}_4(\phi)=\phi,\]
\[\Upsilon_{c_X}(f)=\bigwedge\limits_{x\in X}(c_X(f)(x)\rightarrow f(x)),\,\forall\,f\in L^X,\]
\[\Upsilon_{c_Y}(g)=\bigwedge\limits_{y\in Y}(c_Y(g)(y)\rightarrow g(y)),\,\forall\,g\in L^Y.\]
Now, for all $g\in L^Y$,
\begin{eqnarray*}
\Upsilon_{c_X}(\overleftarrow{\phi}( g))&=&\bigwedge\limits_{x\in X}(c_X(\overleftarrow{\phi}(g))(x)\rightarrow \overleftarrow{\phi}(g)(x))\\
&\geq&\bigwedge\limits_{x\in X}((l\rightarrow \overleftarrow{\phi}(c_Y(g))(x))\rightarrow g(\phi(x)))\\
&=&\bigwedge\limits_{x\in X}((l\rightarrow c_Y(g)(\phi(x)))\rightarrow g(\phi(x)))\\
&\geq& \bigwedge\limits_{x\in X}l\otimes (c_Y(g)(\phi(x))\rightarrow g(\phi(x)))\\
&\geq&  l\otimes \bigwedge\limits_{x\in X}(c_Y(g)(\phi(x))\rightarrow g(\phi(x)))\\
&\geq& l\otimes \bigwedge\limits_{y\in Y}(c_Y(g)(y)\rightarrow g(y))\\
&=& l\otimes \Upsilon_{c_Y}(g).
\end{eqnarray*} 
Thus $\phi:(X,\Upsilon_{c_X})\rightarrow (Y,\Upsilon_{c_Y})$ is a ${\bf FCSS}$-morphism and $F^{-1}_4$ is a functor. Therefore, it can be easily verified that $F_4$ and $F^{-1}_4$ are inverse functors.\\\\
(5) Let $F_5,F^{-1}_5$ be the restriction of $F_4,F^{-1}_4$ to the subcategory ${\bf EFCSS}$ and ${\bf SFCS}$, respectively. It is
clear that $F_5,F^{-1}_5$ are also isomorphism functors.\\\\
It can be verified easily that the diagram of functors in Figure \ref{fig,2} is commutable.
\section{Categorical view of $L$-fuzzy partitions, $T_1$-coalgebras and $(T_2, T_3)$-dialgebras}
In this section, we introduce the concept of coalgebra (dialgebra) corresponding to a space with $L$-fuzzy partition. We begin with the following.\\\\
We study the category ${\bf FPS}$ with the categories of $T_1$-coalgebras and $(T_2, T_3)$-dialgebras. For this purpose, we take an $L$-fuzzy partition under certain conditions, i.e., $L$-fuzzy partition $\mathcal{P}^{'}_X=\{A_{\xi(x)}:x\in X\})$ characterized by index function $\xi:X\rightarrow J=X$ such that $\xi(x)=x,\,\forall\,x\in X.$ The space with this partition is denoted by $(X,\mathcal{P}^{'}_X)$. We shall denote by ${\bf FPS^1}$, the category of spaces with $L$-fuzzy partition $(X,\mathcal{P}^{'}_X)$ and their FP-functions as in Definition \ref{par} ($\phi=\psi$). Now, let $(X,\mathcal{P}^{'}_X)$ be a space with $L$-fuzzy partition $\mathcal{P}^{'}_X$, where $\mathcal{P}^{'}_X=\{A_{x}:x\in X\}$, we denote a function $ F^\uparrow_{
X,\mathcal{P}^{'}_X}:L^X \rightarrow L^X$ such that $F^\uparrow_{X,\mathcal{P}^{'}_X}[f](x) = F^\uparrow_{x}[f],\,\forall\,f\in L^X$.
Next, we have the following from Propositions \ref{1} and \ref{11}.\\\\
Let $(\phi,\phi,\mathcal{W}):(X,\mathcal{P}_X)\rightarrow (Y,\mathcal{P}_Y)$ be a ${\bf FPS^1}$-morphism. Then for all $x\in X,l\in L\setminus 0,g\in L^Y$,
\begin{itemize}
\item[(i)] $ F^\uparrow_{X,\mathcal{P}^{'}_X}[ \overleftarrow{\phi}(g)](x)\otimes {l}\leq  \overleftarrow{\phi}(F^\uparrow_{Y,\mathcal{P}^{'}_Y}[g])(x)$; and
\item[(ii)] $F^\uparrow_{X,\mathcal{P}^{'}_X}[ f](x)\otimes {l}\leq  \overleftarrow{\phi}(F^\uparrow_{Y,\mathcal{P}^{'}_Y}[\overrightarrow{\phi}(f)])(x)$.
\end{itemize}
In the following, we introduce a category of $T_1$-coalgebras.
\begin{pro}\label{T_1} Let $T_1:{\bf SET}\rightarrow {\bf SET}$ be a function such that for all $X\in|{\bf SET}|$, $T_1(X)=L^{L^X}$ and for all {\bf SET}-morphism $\phi:X\rightarrow Y$, there exists a function $T_1(\phi):L^{L^X}\rightarrow L^{L^Y}$ such that for all $\lambda\in L^{L^X}, \,T_1(\phi)(\lambda)=\overleftarrow{\overleftarrow{\phi}}(\lambda)$, where $\overleftarrow{\overleftarrow{\phi}}$ is backward operator from $L^{L^X}$ to $L^{L^Y}$. Then $T_1$ is functor.
\end{pro}
\textbf{Proof:} (i) Let $\lambda\in L^{L^X},f\in L^Z$ and $\phi_1:X\rightarrow Y,\phi_2:Y\rightarrow Z$ be {\bf SET}-morphisms. Then
\begin{eqnarray*}
(T_1(\phi_2)\circ T_1(\phi_1))(\lambda)(f)&=&T_1(\phi_2)(T_1(\phi_1)(\lambda))(f)\\
&=&\overleftarrow{\overleftarrow{\phi_2}}(T_1(\phi_1)(\lambda))(f)\\
&=&T_1(\phi_1)(\lambda)(\overleftarrow{\phi_2}(f))\\
&=&\overleftarrow{\overleftarrow{\phi_1}}(\lambda)(\overleftarrow{\phi_2}(f))\\
&=&\lambda(\overleftarrow{\phi_1}(\overleftarrow{\phi_2}(f)))\\
&=&\lambda((\overleftarrow{\phi_1}\circ\overleftarrow{\phi_2})(f))\\
&=&\lambda(\overleftarrow{{\phi_2}\circ{\phi_1}}(f))\\
&=&\overleftarrow{\overleftarrow{{\phi_2}\circ{\phi_1}}}(\lambda)(f)\\
&=&T_1({\phi_2}\circ{\phi_1})(\lambda)(f).
\end{eqnarray*}
Thus $T_1(\phi_2)\circ T_1(\phi_1)=T_1({\phi_2}\circ{\phi_1})$.\\
(ii) Let $X \in |{\bf SET}|, \lambda \in L^{L^X}, f\in L^X$. Then $T_1(id_X)(\lambda)(f)=\overleftarrow{\overleftarrow{id_X}}(\lambda)(f)=\lambda(f).$ Thus $T_1(id_X)(\lambda)=\lambda$, which implies that $T_1(id_X) = id_{T_1(X)}$. Hence $T_1$ is a functor.
\begin{def1} For a functor $T_1:{\bf SET}\rightarrow {\bf SET}$, $T_1$-{\bf coalgebra} is a pair $\mathcal{X}_C=(X,\alpha_X)$, where $X\in|{\bf SET}|$ and $\alpha_X:X\rightarrow T_1(X)$, i.e., $\alpha_X:X\rightarrow L^{L^X}$ is a structure function of $\mathcal{X}_C$.
\end{def1}
\begin{def1}
Let $\mathcal{X}_C=(X,\alpha_X)$ and $\mathcal{X}^{'}_C=(Y,\alpha_Y)$ be $T_1$-{coalgebras}. Then  $\phi:\mathcal{X}_C \rightarrow \mathcal{X}^{'}_C$ is a {\bf homomorphism} if $\phi:X\rightarrow Y$ is a function such that $(T_1(\phi)\circ \alpha_X)\leq \alpha_Y\circ \phi$.
\end{def1} 
\begin{pro}\label{T_1c} $T_1$-coalgebras alongwith their homomorphisms form a category.
\end{pro}
\textbf{Proof:} We only need to show that composition of homomorphisms is also a homomorphism. For which, let $\phi_1:(X,\alpha_X)\rightarrow (Y,\alpha_Y)$ and $\phi_2:(Y,\alpha_Y)\rightarrow (Z,\alpha_Z)$ be homomorphisms, i.e., $\phi_1:X\rightarrow Y$ and $\phi_2:Y\rightarrow Z$ are functions  such that $(T_1(\phi_1)\circ\alpha_X )\leq \alpha_Y\circ \phi_1,(T_1(\phi_2)\circ\alpha_Y )\leq \alpha_Z\circ \phi_2$. Then
\begin{eqnarray*}
(T_1(\phi_2\circ\phi_1)\circ \alpha_X)(x)&=&((T_1(\phi_2)\circ T_1(\phi_1))\circ\alpha_X)(x)\\
&=&(T_1(\phi_2)\circ (T_1(\phi_1)\circ\alpha_X))(x)\\
&\leq& (T_1(\phi_2)\circ (\alpha_Y\circ \phi_1))(x)\\
&=& T_1(\phi_2) ((\alpha_Y\circ \phi_1)(x))\\
&=& T_1(\phi_2) (\alpha_Y( \phi_1(x)))\\
&=& (T_1(\phi_2)\circ\alpha_Y)( \phi_1(x))\\
&\leq&(\alpha_Z\circ \phi_2)(\phi_1(x))\\
&=&((\alpha_Z \circ\phi_2)\circ\phi_1)(x)\\
&=& (\alpha_Z\circ (\phi_2\circ\phi_1))(x).
\end{eqnarray*}
Therefore $T_1(\phi_2\circ\phi_1)\circ \alpha_X\leq \alpha_Z\circ (\phi_2\circ\phi_1)$. Thus $\phi_2\circ\phi_1:(X,\alpha_X)\rightarrow (Z,\alpha_Z)$ is a homomorphism.\\\\
We shall denote by ${\bf COA}$, the category of $T_1$-coalgebras and their homomorphisms. Next, we have the following.\\\\
Let $(X,\mathcal{P}^{'}_X)\in |{\bf FPS^1}|$. Then the upper $F$-transform $F^\uparrow_{X,\mathcal{P}^{'}_X}:L^X\rightarrow L^X$ of $f$ may be interpreted as $T_1$-coalgebra structure function $\alpha_{\mathcal{P}^{'}_X}:X\rightarrow L^{L^X}$ such that $\alpha_{\mathcal{P}^{'}_X}(x)(f)=F_{X,\mathcal{P}^{'}_X}^\uparrow[f](x),\,\forall\,x\in X,f\in L^X$. Thus  $(X,\mathcal{P}^{'}_X)$ can be viewed as a $T_1$-coalgebra ${\mathcal{P}^{'}_{X_C}}=(X,\alpha_{\mathcal{P}^{'}_X})$.
\begin{pro}\label{T_1fp}
If $(\phi,\phi,\mathcal{W}):(X,\mathcal{P}^{'}_X)\rightarrow(Y,\mathcal{P}^{'}_Y)$ is a ${\bf FPS^1}$-morphism. Then $\phi:(X,\alpha_{\mathcal{P}^{'}_X})\rightarrow (Y,\alpha_{\mathcal{P}^{'}_Y})$ is a ${\bf COA}$-morphism.
\end{pro}
\textbf{Proof:} Let $x\in X,f\in L^Y$ and $(\phi,\phi,\mathcal{W}):(X,\mathcal{P}^{'}_X)\rightarrow (Y,\mathcal{P}^{'}_Y)$ be a ${\bf FPS^1}$-morphism. Then
\begin{eqnarray*}
(T_1(\phi)\circ\alpha_{\mathcal{P}^{'}_X})(x)(f)&=&T_1(\phi)(\alpha_{\mathcal{P}^{'}_X}(x))(f)\\
&=&\overleftarrow{\overleftarrow{\phi}}(\alpha_{\mathcal{P}^{'}_X}(x))(f)\\
&=&\alpha_{\mathcal{P}^{'}_X}(x)(\overleftarrow{\phi}(f))\\
&=&F^\uparrow_{X,\mathcal{P}^{'}_X}[\overleftarrow{\phi}(f)](x)\\
&=&1\otimes F^\uparrow_{X,\mathcal{P}^{'}_X}[\overleftarrow{\phi}(f)](x)\\
&\leq& \overleftarrow{\phi}(F^\uparrow_{Y,\mathcal{P}^{'}_Y}[f])(x)\\
&=& F^\uparrow_{Y,\mathcal{P}^{'}_Y}[f](\phi(x))\\
&=&\alpha_{\mathcal{P}^{'}_Y}(\phi(x))(f)\\
&=&( \alpha_{\mathcal{P}^{'}_Y}\circ\phi)(x)(f).
\end{eqnarray*}
Therefore $(T_1(\phi)\circ\alpha_{\mathcal{P}^{'}_X})\leq \alpha_{\mathcal{P}^{'}_Y}\circ\phi$. Thus $\phi:(X,\alpha_{\mathcal{P}^{'}_X})\rightarrow (Y,\alpha_{\mathcal{P}^{'}_Y})$ is a ${\bf COA}$-morphism.
\begin{pro}
Let $F: {\bf FPS^1}\rightarrow {\bf COA}$ be a function such that for all $(X,\mathcal{P}^{'}_X)\in|{\bf FPS^1}|,\,F(X,\mathcal{P}^{'}_X)=(X,\alpha_{\mathcal{P}^{'}_X})$ and for every ${\bf FPS^1}$-morphism $(\phi,\phi,\mathcal{W}):(X,\mathcal{P}^{'}_X)\rightarrow (Y,\mathcal{P}^{'}_Y)$, $F(\phi,\phi,\mathcal{W})
:(X,\alpha_{\mathcal{P}^{'}_X})\rightarrow (Y,\alpha_{\mathcal{P}^{'}_Y})$ such that $F(\phi,\phi,\mathcal{W})=\phi$. Then $F$ is a functor.
\end{pro}
\begin{pro}\label{T_2} Let $T_2 : {\bf SET} \rightarrow {\bf SET}$ be a function such that for all $X \in|{\bf SET}|, T_2(X) = X \times L^X$ and for every ${\bf SET}$-morphism $\phi :X \rightarrow Y$, $T_2(\phi):X\times L^X\rightarrow Y\times L^Y$ such that $T_2(\phi)=(\phi,\overrightarrow{\phi})$. Then $T_2$ is a functor.
\end{pro}
\textbf{Proof:} (i) Let $x\in X,f\in L^X$ and $\phi_1:X\rightarrow Y,\phi_2:Y\rightarrow Z$ be ${\bf SET}$-morphisms. Then
\begin{eqnarray*}
(T_2(\phi_2)\circ T_2(\phi_1))(x,f)&=&T_2(\phi_2)(T_2(\phi_1)(x,f))\\
&=&T_2(\phi_2)((\phi_1,\overrightarrow{\phi_1})(x,f))\\
&=&(\phi_2,\overrightarrow{\phi_2})((\phi_1,\overrightarrow{\phi_1})(x,f))\\
&=&(\phi_2,\overrightarrow{\phi_2})(\phi_1(x),\overrightarrow{\phi_1}(f))\\
&=&(\phi_2(\phi_1(x)),\overrightarrow{\phi_2}(\overrightarrow{\phi_1}(f)))\\
&=&((\phi_2\circ\phi_1)(x),(\overrightarrow{\phi_2}\circ\overrightarrow{\phi_1}) (f))\\
&=&((\phi_2\circ\phi_1)(x),\overrightarrow{\phi_2\circ\phi_1} (f))\\
&=&(\phi_2\circ\phi_1,\overrightarrow{\phi_2\circ\phi_1})(x,f)\\
&=& T_2(\phi_2\circ\phi_1)(x,f)
\end{eqnarray*}
Thus $T_2(\phi_2)\circ T_2(\phi_1)=T_2(\phi_2\circ\phi_1)$.\\
(ii) Let $x\in X,f\in L^X$. Then $T_2(id_X)(x,f)=(id_X,\overrightarrow{id_{X}})(x,f)=(id_X(x),\overrightarrow{id_{X}}(f))$ $=(x,f)$. Therefore $T_2(id_X)=id_{T_2(X)}$. Thus $T_2$ is a functor.
\begin{pro} \label{T_3}
Let $T_3 : {\bf SET }\rightarrow {\bf SET}$ be a function such that for all $X \in |{\bf SET}|$, $T_3(X) = L$ and for every ${\bf SET}$-morphism $\phi :X\rightarrow Y$, $T_3(\phi):L\rightarrow L$ be a function such that for all $T_3(\phi)=id_L$. Then $T_3$ is a functor.
\end{pro}
\textbf{Proof:} (i) Let $a\in L$ and $\phi_1:X\rightarrow Y,\phi_2:Y\rightarrow Z$ be ${\bf SET}$-morphisms. Then
\begin{eqnarray*}
(T_3(\phi_2)\circ T_3(\phi_1))(a)&=&T_3(\phi_2)(T_3(\phi_1)(a))\\
&=&id_L(T_3(\phi_1)(a))\\
&=&id_L(id_L(a))\\
&=&id_L(a)\\
&=&T_3(\phi_2\circ\phi_1)(a).
\end{eqnarray*}
Thus $T_3(\phi_2)\circ T_3(\phi_1)=T_3(\phi_2\circ\phi_1)$.\\
(ii) Let $a\in L$. Then $T_3(id_X)(a)=id_L(a)=a$. Therefore $T_3(id_X)=id_{T_3(X)}$. Thus $T_3$ is a functor.
\begin{def1} For functors $T_2, T_3 : {\bf SET}\rightarrow {\bf SET}$, a $(T_2, T_3)$-{\bf dialgebra} is a pair $\mathcal{X}_D = (X, \beta_X)$, where $X \in |{\bf SET}|$ and $\beta_X : T_2(X) \rightarrow T_3(X)$, i.e., $\beta_X: X\times L^X\rightarrow L$ is a structure function of $\mathcal{X}_D$.
\end{def1}
In the following, we introduce the category of $(T_2,T_3)$-dialgebras.
\begin{def1} Let $\mathcal{X}_D = (X, \beta_X)$ and $\mathcal{X}^{'}_D=(Y,\beta_Y)$ be the $(T_2, T_3)$-dialgebras. Then $\phi:\mathcal{X}_D\rightarrow\mathcal{X}^{'}_D$ is a {\bf homomorphism} if $\phi : X \rightarrow Y$ is a function such that $T_3(\phi)\circ \beta_X\leq \beta_Y\circ T_2(\phi)$.
\end{def1}
\begin{pro}\label{T_2c}
$(T_2, T_3$)-dialgebras alongwith their homomorphisms form a category.
\end{pro}
\textbf{Proof:} Similar to that of Proposition \ref{T_1c}.\\\\
We shall denote by ${\bf DIA}$, the category of $(T_2,T_3)$-dialgebras and their homomorphisms.\\
Let $(X,\mathcal{P}^{'}_X)\in |{\bf FPS^1}|$. Then the upper $F$-transform $F^\uparrow_{X,\mathcal{P}^{'}_X}:L^X\rightarrow L^X$ of $f$ may be interpreted as $(T_2,T_3)$-dialgebra structure function $\beta_{\mathcal{P}^{'}_X}:X\times L^X\rightarrow L$ such that $\beta_{\mathcal{P}^{'}_X}(x,f)=F_{X,\mathcal{P}^{'}_X}^\uparrow[f](x),\,\forall\,x\in X,f\in L^X$. Thus  $(X,\mathcal{P}^{'}_X)$ can be viewed as a $(T_2,T_3)$-dialgebra $\mathcal{P}^{'}_{X_D}=(X,\beta_{\mathcal{P}^{'}_X})$.
\begin{pro}\label{T_2fp}
If $(\phi,\phi,\mathcal{W}):(X,\mathcal{P}^{'}_X)\rightarrow(Y,\mathcal{P}^{'}_Y)$ is a ${\bf FPS^1}$-morphism. Then $\phi:(X,\beta_{\mathcal{P}^{'}_X})\rightarrow (Y,\beta_{\mathcal{P}^{'}_Y})$ is a ${\bf DIA}$-morphism.
\end{pro}
\textbf{Proof:} Let $x\in X,f\in L^X$ and $(\phi,\phi,\mathcal{W}):(X,\mathcal{P}^{'}_X)\rightarrow (Y,\mathcal{P}^{'}_Y)$ be a ${\bf FPS^1}$-morphism. Then
\begin{eqnarray*} 
(\beta_{\mathcal{P}^{'}_Y}\circ T_2(\phi))(x,f)&=&\beta_{\mathcal{P}^{'}_Y}(T_2(\phi)(x,f))\\
&=&\beta_{\mathcal{P}^{'}_Y}((\phi,\overrightarrow{\phi})(x,f))\\
&=&\beta_{\mathcal{P}^{'}_Y}(\phi(x),\overrightarrow{\phi}(f))\\
&=&F^\uparrow_{Y,\mathcal{P}^{'}_Y}[\overrightarrow{\phi}(f)](\phi(x))\\
&=&\overleftarrow{\phi}(F^\uparrow_{Y,\mathcal{P}^{'}_Y}[\overrightarrow{\phi}(f)])(x)\\
&=&1\rightarrow\overleftarrow{\phi}(F^\uparrow_{Y,\mathcal{P}^{'}_Y}[\overrightarrow{\phi}(f)])(x)\\
&\geq&F^\uparrow_{X,\mathcal{P}^{'}_X}[f](x)\\
&=&\beta_{\mathcal{P}^{'}_X}(x,f)\\
&=&id_L(\beta_{\mathcal{P}^{'}_X}(x,f))\\
&=& (T_3(\phi)\circ\beta_{\mathcal{P}^{'}_X})(x,f).
\end{eqnarray*}
Therefore $T_3(\phi)\circ\beta_{\mathcal{P}^{'}_X}\leq\beta_{\mathcal{P}^{'}_Y}\circ T_2(\phi)$. Thus $\phi:(X,\beta_{\mathcal{P}^{'}_X})\rightarrow (Y,\beta_{\mathcal{P}^{'}_Y})$ is a ${\bf DIA}$-morphism.
\begin{pro}
Let $F^{'}: {\bf FPS^1}\rightarrow {\bf DIA}$ be a function such that for all $(X,\mathcal{P}^{'}_X)\in|{\bf FPS^1}|,\,F(X,\mathcal{P}^{'}_X)=(X,\beta_{\mathcal{P}^{'}_X})$ and for every ${\bf FPS^1}$-morphism $(\phi,\phi,\mathcal{W}):(X,\mathcal{P}^{'}_X)\rightarrow (Y,\mathcal{P}^{'}_Y)$, $F(\phi,\phi,\mathcal{W}):(X,\beta_{\mathcal{P}^{'}_X})\rightarrow (Y,\beta_{\mathcal{P}^{'}_Y})$ such that $F(^{'}\phi,\phi,\mathcal{W})=\phi$. Then $F^{'}$ is a functor.
\end{pro}
In the following, we establish an isomorphism between the category of $T_1$-coalgebras and the category of $(T_2, T_3$)-dialgebras. Let $\mathcal{X}_C = (X, \alpha_X)$ be a $T_1$-coalgebra and $\mathcal{X}_D = (X, \beta_X)$ be a $(T_2, T_3$)-dialgebra. Then we construct a $(T_2, T_3$)-dialgebra and $T_1$-coalgebra corresponding to $\mathcal{X}_C$ and $\mathcal{X}_D$, respectively, denotded by $(X,\beta_{\alpha_X})$ and $(X,\alpha_{\beta_X})$, as follows:
\begin{itemize}
\item[(i)]  ${\mathcal{X}_{C}}_{D} = (X, \beta_{\alpha_X})$, where  $\beta_{\alpha_X}:X\times L^X\rightarrow L$ is a structure function of ${\mathcal{X}_{C}}_{D}$ such that $\beta_{\alpha_X}(x, f) = \alpha_X(x)(f),\,\forall\,x\in X,f\in L^X$; and
\item[(ii)] ${\mathcal{X}_{D}}_{C} = (X, \alpha_{\beta_X})$, where  $\alpha_{\beta_X}:X\rightarrow L^{L^X}$ is a structure function of ${\mathcal{X}_{D}}_{C}$ such that $\alpha_{\beta_X}(x)( f) = \beta_X(x,f),\,\forall\,x\in X,f\in L^X$.
\end{itemize}
To show isomorphism, we need the following result from \cite{ying}.
\begin{pro} Let $\phi:X\rightarrow Y$ be a function. Then 
\begin{itemize}
\item[(i)] for all $f\in L^X,\,(\overleftarrow{\phi}\circ\overrightarrow{\phi})(f)=f$, if $\overrightarrow{\phi}:L^X\rightarrow L^Y$ is injective; and
\item[(ii)] for all $g\in L^Y,\,(\overrightarrow{\phi}\circ\overleftarrow{\phi})(g)=g$, if $\overrightarrow{\phi}:L^X\rightarrow L^Y$ is surjective.
    \end{itemize}
    \end{pro}
\begin{pro}
 Every ${\bf COA}$-morphism is a ${\bf DIA}$-morphism, provided that $\overrightarrow{\phi}$ is injective.
\end{pro}
\textbf{Proof:} Let $\mathcal{X}_C=(X,\alpha_X),\mathcal{X}^{'}_C=(Y,\alpha_{Y})\in|{\bf COA}|$ and $\phi:\mathcal{X}_C\rightarrow \mathcal{X}^{'}_C$ be a ${\bf COA}$-morphism. Then for all $x\in X,f\in L^X$
\begin{eqnarray*}
(\beta_{\alpha_Y}\circ T_2(\phi))(x,f)&=&\beta_{\alpha_Y}( T_2(\phi)(x,f))\\
&=&\beta_{\alpha_Y}((\phi,\overrightarrow{\phi})(x,f))\\
&=&\beta_{\alpha_Y}(\phi(x),\overrightarrow{\phi}(f))\\
&=&\alpha_Y(\phi(x))(\overrightarrow{\phi}(f))\\
&=&(\alpha_Y\circ\phi)(x)(\overrightarrow{\phi}(f))\\
&\geq& (T_1(\phi)\circ \alpha_X)(x)(\overrightarrow{\phi}(f))\\
&=& T_1(\phi)(\alpha_X(x))(\overrightarrow{\phi}(f))\\
&=& \overleftarrow{\overleftarrow{\phi}}(\alpha_X(x))(\overrightarrow{\phi}(f))\\
&=& \alpha_X(x)(\overleftarrow{\phi}(\overrightarrow{\phi}(f)))
\\
&=& \alpha_X(x)(f)\\
&=& \beta_{\alpha_X}(x,f)\\
&=& id_L(\beta_{\alpha_X}(x,f))\\
&=&(T_3(\phi)\circ\beta_{\alpha_X})(x,f).
\end{eqnarray*}
Thus $\beta_{\alpha_Y}\circ T_2(\phi)\geq (T_3(\phi)\circ\beta_{\alpha_X})$. Therefore $\phi:{\mathcal{X}_{C}}_{D}\rightarrow {\mathcal{X}^{'}_{C}}_{D}$ is a ${\bf DIA}$-morphism.
\begin{pro}
Every ${\bf DIA}$-morphism is a ${\bf COA}$-morphism, provided that $\overrightarrow{\phi}$ is surjective.
\end{pro}
\textbf{Proof:} Let $\mathcal{X}_D=(X,\beta_X),\mathcal{X}^{'}_D=(Y,\beta_{Y})\in|{\bf DIA}|$ and $\phi:\mathcal{X}_D\rightarrow \mathcal{X}^{'}_D$ be a ${\bf DIA}$-morphism. Then for all $x\in X,f\in L^X$
\begin{eqnarray*}
(T_1(\phi)\circ \alpha_{\beta_X})(x)(f)&=&T_1(\phi) (\alpha_{\beta_X}(x))(f)\\
&=&\overleftarrow{\overleftarrow{\phi}}(\alpha_{\beta_X}(x))(f)
\end{eqnarray*}
\begin{eqnarray*}
&=&(\alpha_{\beta_X}(x))(\overleftarrow{\phi}(f))\\
&=&{\beta_X}(x,\overleftarrow{\phi}(f))\\
&=&id_L({\beta_X}(x,\overleftarrow{\phi}(f)))\\
&=& (T_3(\phi)\circ \beta_X)(x,\overleftarrow{\phi}(f))\\
&\leq& (\beta_Y\circ T_2(\phi))(x,\overleftarrow{\phi}(f))\\
&=&\beta_Y(T_2(\phi)(x,\overleftarrow{\phi}(f)))\\
&=&\beta_Y((\phi,\overrightarrow{\phi})(x,\overleftarrow{\phi}(f)))\\
&=&\beta_Y(\phi(x),\overrightarrow{\phi}(\overleftarrow{\phi}(f)))\\
&=&\beta_Y(\phi(x),f)\\
&=& \alpha_{\beta_Y}(\phi(x))(f)\\
&=& (\alpha_{\beta_Y}\circ\phi)(x)(f).
\end{eqnarray*}
Thus $T_1(\phi)\circ \alpha_{\beta_X}\leq \alpha_{\beta_Y}\circ\phi$. Therefore $\phi:{\mathcal{X}_D}_C\rightarrow {\mathcal{X}^{'}_D}_C$ is a ${\bf COA}$-morphism.
\begin{pro}\label{T}
Let $T:{\bf COA}\rightarrow {\bf DIA}$ be a function such that for every $\mathcal{X}_C\in|{\bf COA}|,$ $T(\mathcal{X}_C)={\mathcal{X}_{C}}_{D},$ for every ${\bf COA}$-morphism $\phi:\mathcal{X}_C\rightarrow \mathcal{X}^{'}_C,T(\phi):{\mathcal{X}_{C}}_{D}\rightarrow {\mathcal{X}^{'}_{C}}_{D}$ be a function such that $T(\phi)=\phi$. Then $T$ is a functor.
\end{pro}
\textbf{Proof:} (i) Let $\phi:\mathcal{X}_C\rightarrow \mathcal{X}^{'}_C,\phi^{'}:\mathcal{X}^{'}_C\rightarrow \mathcal{X}^{''}_C$ be ${\bf COA}$-morphisms, where $\mathcal{X}_C=(X,\alpha_X),\mathcal{X}^{'}_C=(X^{'},\alpha_{X^{'}})$ and $\mathcal{X}^{''}_C=(X^{''},\alpha_{X^{''}})$. Then
\[T(\phi^{'}\circ\phi)=\phi^{'}\circ\phi=T(\phi^{'})\circ T(\phi).\] 
Thus $T(\phi^{'}\circ\phi)=T(\phi^{'})\circ T(\phi)$.\\\\
(ii) Let $\mathcal{X}_C=(X,\alpha_X)\in|{\bf COA}|$. Then $T(id_{\mathcal{X}_C})=T(id_X)=id_{X}$. Thus $T(id_{\mathcal{X}_C})=id_{T(\mathcal{X}_C)}$
\begin{pro}
Let $T^{'}:{\bf DIA}\rightarrow {\bf COA}$ be a function such that for every $\mathcal{X}_D\in|{\bf DIA}|,T^{'}(\mathcal{X}_D)={\mathcal{X}_{D}}_{C}$, for every ${\bf DIA}$-morphism $\phi:\mathcal{X}_D\rightarrow \mathcal{X}^{'}_D,T^{'}(\phi):{\mathcal{X}_{D}}_C\rightarrow {\mathcal{X}^{'}_{D}}_{C}$ be a function such that $T^{'}(\phi)=\phi$. Then $T^{'}$ is a functor.
\end{pro}
\textbf{Proof:} Similar to that of Proposition \ref{T}.
\begin{pro}
The category ${\bf COA}$ is isomorphic to the category ${\bf DIA}$.
\end{pro}
\textbf{Proof:} Let $T:{\bf COA}\rightarrow {\bf DIA},T^{'}:{\bf DIA}\rightarrow {\bf COA}$ be functors. Then $T^{'}\circ T:{\bf COA}\rightarrow {\bf COA}$ is a functor such that for every $\mathcal{X}_C\in|{\bf COA}|$ and for every ${\bf COA}$-morphism $\phi:\mathcal{X}_C\rightarrow \mathcal{X}^{'}_C$,
\[(T^{'}\circ T)(\mathcal{X}_C)=T^{'}( T(\mathcal{X}_C))=T^{'}({\mathcal{X}_C}_D)={{\mathcal{X}_C}_D}_D={\mathcal{X}_C}.\]
Thus $(T^{'}\circ T)(\mathcal{X}_C)=\mathcal{X}_C$ and 
\[(T^{'}\circ T)(\phi)=T^{'}( T(\phi))=T^{'}(\phi)=\phi.\] 
Thus  $(T^{'}\circ T)(\phi)=\phi$. Therefore $T^{'}\circ T=id_{\bf COA}$. Similarly, we can show that $T\circ T^{'}:{\bf DIA}\rightarrow {\bf DIA}$ is a functor such that $T\circ T^{'}=id_{\bf DIA}$. Hence the category ${\bf COA}$ is isomorphic to the category ${\bf DIA}$.
Next, we have the following.
\begin{pro} \label{com}
Let $F:{\bf FPS}^1\rightarrow {\bf COA}$, $F^{'}:{\bf FPS}^1\rightarrow {\bf DIA}$ and $T:{\bf COA}\rightarrow {\bf DIA}$ be functors. Then the diagram in Figure \ref{fig,6} is commutable. 
\end{pro}
\textbf{Proof:} Let $\mathcal{P}^{'}_X\in |{\bf FPS}^1|$ and $(\phi,\phi,\mathcal{W}):(X,\mathcal{P}^{'}_X)\rightarrow (X^{'},\mathcal{P}^{'}_{X)^{'}})$ be a ${\bf FPS}^1$-morphism. Then
\[(T\circ F)(X,\mathcal{P}^{'}_X)=T(F(X,\mathcal{P}^{'}_X)) =T(X,\alpha_{\mathcal{P}^{'}_X})=(X,\beta_{\alpha_{\mathcal{P}^{'}_X}})=F^{'}(X,\mathcal{P}^{'}_X).\] 
Thus $(T\circ F)(X,\mathcal{P}^{'}_X)=F^{'}(X,\mathcal{P}^{'}_X)$, and 
\[(T\circ F)(\phi,\phi,\mathcal{W})=T(F(\phi,\phi,\mathcal{W})) =T(\phi)=\phi=F^{'}(\phi,\phi,\mathcal{W}).\]
Thus $(T\circ F)(\phi,\phi,\mathcal{W})=F^{'}(\phi,\phi,\mathcal{W})$. Hence $T\circ F=F^{'}$, i.e., diagram in Figure \ref{fig,6} is commutable.
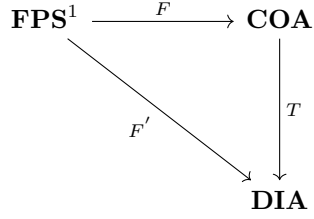
\begin{figure}
\[
 \begin{tikzcd}[row sep=12ex, column sep=12ex] 
   {\bf FPS}^{1}\arrow{r}{F} \arrow[swap]{dr}{F^{'}} & {\bf COA}\arrow{d}{T} \\
     & {\bf DIA}
  \end{tikzcd}
  \]
 \caption{Diagram for Proposition \ref{com}}
\label{fig,6}
  \end{figure}\\
 Finally, we show that there exists a pair of functors $(T,T^{'})$ that has adjoint property.
\begin{pro}\label{adj1}
Let $T:{\bf COA}\rightarrow {\bf DIA}$ and $T^{'}:{\bf DIA}\rightarrow {\bf COA}$ be functors. Then $T$ is a left adjoint to $T^{'}$ and $T^{'}$ is a right adjoint to $T$.
\end{pro}
\textbf{Proof:} To show this result, we have to show that there exists a natural transformation $\Psi:id_{\bf COA}\rightarrow T^{'}\circ T$ such that for every $\mathcal{X}_C\in |{\bf COA}|$ and ${\bf COA}$-morphism $\phi:\mathcal{X}_C\rightarrow T^{'}(\mathcal{X}_D)$, there exists a unique ${\bf DIA}$-morphism $\rho:T(\mathcal{X}_C)\rightarrow \mathcal{X}_D$ such that the diagram in Figure \ref{fig,8} commutes. For which, let $\mathcal{X}_C=(X,\alpha_X)\in |{\bf COA}|,\mathcal{X}_D=(X,\beta_X)\in|{\bf DIA}|$ and $\Psi_{\mathcal{X}_C}:(X,\alpha_X)\rightarrow (T^{'}\circ T)(X,{\alpha_X})$ be a function such that $\Psi_{\mathcal{X}_C}=id_X$. It is easy to check that $\Psi$ is a natural transformation. Also, let $\phi:(X,\alpha_X)\rightarrow T^{'}(X,\beta_X)$ be a ${\bf COA}$-morphism. Now, we define a ${\bf DIA}$-morphism $\rho:T(\mathcal{X}_C)\rightarrow \mathcal{X}_D$ such that $\rho=\phi$. The diagram in Figure \ref{fig,9} commutes, i.e., $T^{'}(\rho)\circ id_X=\rho\circ id_X=\rho=\phi$.
Thus $T^{'}(\rho)\circ id_X=\phi$. The uniqueness of $\rho$ is trivial. Hence $T$ is a left adjoint of $T^{'}$ and $T^{'}$ is a right adjoint of $T$.
\begin{figure}
\[
 \begin{tikzcd}[row sep=12ex, column sep=12ex] 
   \mathcal{X}_C\arrow{r}{\Psi_{\mathcal{X}_C}} \arrow[swap]{dr}{\phi} & T^{'}(T(\mathcal{X}_C))\arrow{d}{T^{'}(\rho)} &T(\mathcal{X}_C)\arrow{d}{\rho}\\
     & T^{'}(\mathcal{X}_D)&\mathcal{X}_D
  \end{tikzcd}
  \]
 \caption{Diagram for Proposition \ref{adj1}}
\label{fig,8}
  \end{figure}
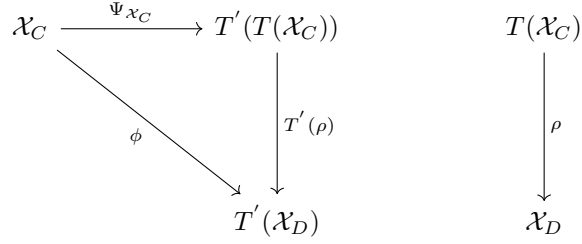
  \begin{figure}
    \[
  \begin{tikzcd}[row sep=12ex, column sep=10ex]  
  (X,\alpha_X)\arrow{r}{id_X} \arrow[swap]{dr}{\phi} & T^{'}(T(X,\alpha_X))\arrow{d}{T^{'}(\rho)} &T(X,\alpha_X)\arrow{d}{\rho}\\
     & T^{'}(X,\beta_X)&(X,\beta_X)
    \end{tikzcd}
\]
\caption{Diagram is equivalent to the diagram in Figure \ref{fig,8}}
\label{fig,9}
  \end{figure}
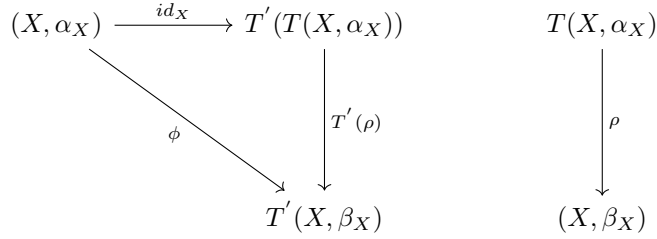
\section{Concluding remarks}
{This contribution is towards using different theories to enrich the theory of $F$-transform. {Specifically, we have demonstrated that the direct upper $F$-transform determines an $L$-fuzzy closure system uniquely and introduced $L$-fuzzy partition spaces, $L$-fuzzy closure system spaces in a categorical framework}. Moreover, we have shown a functorial relationship among the introduced categories ${\bf {FPS}},{\bf FAS},\textbf{FCSS}$ and $\textbf{FCS}$. Interestingly, the categories of spaces with $L$-fuzzy partitions, $L$-fuzzy approximation spaces, $L$-fuzzy system spaces and $L$-fuzzy closure spaces studied in \cite{bel,fan,lai,mo} turn out to be  particular cases of the categories ${\bf {FPS}},{\bf FAS},\textbf{FCSS}$ and $\textbf{FCS}$, respectively. Further, we have studied the concept of coalgebra (dialgebra) {corresponding to a direct upper $F$-transform} and shown the existence of functors between the category ${\bf FPS^1}$ and the categories ${\bf COA}$ and ${\bf DIA}$, respectively. Finally, it is shown here that there exists an isomorphism between the categories ${\bf COA}$ and ${\bf DIA}$ under certain conditions alongwith adjointness property.} 


\end{document}